\documentclass[11pt, a4paper]{article}

\usepackage{amsthm, amssymb, amsmath, textcomp, mathrsfs}
\usepackage[latin1,ansinew]{inputenc}
\usepackage[pdftex]{graphics}
\usepackage{graphicx}
\usepackage{longtable}
\usepackage{scrpage2}
\usepackage{hyphenat}
\usepackage{paralist}
\usepackage{color}

\makeatletter
\renewcommand{\section}{%
  \@startsection {section}{1}{\z@}%
                 {-3.5ex plus -1ex minus -.2ex}%
                 {2.3ex plus.2ex}%
                 {\normalfont\Large\bfseries}}

\renewcommand{\subsection}{%
  \@startsection {subsection}{1}{\z@}%
                 {-3.5ex plus -1ex minus -.2ex}%
                 {2.3ex plus.2ex}%
                 {\normalfont\large\bfseries}}
                 
\renewcommand{\subsubsection}{%
  \@startsection {subsubsection}{1}{\z@}%
                 {-2.5ex plus -1ex minus -.2ex}%
                 {1.3ex plus.2ex}%
                 {\normalfont\bfseries}}
\makeatother

\setheadsepline{.4pt}
\clearscrheadfoot
\lehead{\parbox[b]{0.8cm}{\pagemark}\textsc{Chapter \thechapter}}
\rohead{\textsc{\rightmark}\parbox[b]{0.8cm}{\hfill \pagemark}}

	\newcommand*{\savesymbol}[1]{%
	  \expandafter\let\csname orig#1\expandafter\endcsname\csname#1\endcsname
	  \expandafter\let\csname #1\endcsname\relax
	}
	
	\newcommand*{\restoresymbol}[2]{%
	  \expandafter\global\expandafter\let\csname#1#2\expandafter\endcsname%
	    \csname#2\endcsname
	  \expandafter\global\expandafter\let\csname#2\expandafter\endcsname%
	    \csname orig#2\endcsname
	}
	
	\newcommand*{\renamerobustsymbol}[2]{%
	  \expandafter\let\expandafter\origrealcommand
	    \csname #2\space\endcsname
	  \expandafter\global\expandafter\let\csname#1#2\endcsname=\origrealcommand
	}
	
	\makeatletter
	\def\ifnotsavedsym@helper#1#2!{\expandafter\ifx\csname orig#2\endcsname\relax}
	\newcommand*{\ifnotsavedsym}[1]{%
	  \expandafter\ifnotsavedsym@helper\string#1!%
	}
	\makeatother
	
	\newif\ifloadpackages
	\loadpackagestrue
	
	\makeatletter
	\newcommand{\missingpkgs}{}
	\newcommand{\foundpkgs}{}
	\newcommand{\if@sty@file@exists@star}[3]{%
	  \ifloadpackages
	    \IfFileExists{#1.sty}{#2}{#3}%
	  \else
	    #3%
	  \fi
	}
	\newcommand{\if@sty@file@exists}[3]{%
	  \ifloadpackages
	    \IfFileExists{#1.sty}%
	                 {#2\@cons\foundpkgs{{#1}}}%
	                 {#3\completefalse\@cons\missingpkgs{{#1}}}%
	  \else
	    #3\completefalse\@cons\missingpkgs{{#1}}%
	  \fi
	}
	\newcommand{\IfStyFileExists}{%
	  \@ifstar{\if@sty@file@exists@star}{\if@sty@file@exists}%
	}
	\makeatother

	\allowdisplaybreaks
	
	\newcommand{\dStr}[1]{{\rm \textbf{#1}}}	

	\usepackage{stmaryrd}
	
	\newif\ifMTOOLS
	
	\IfStyFileExists{mathtools}
  {\MTOOLStrue
   \savesymbol{xleftrightarrow} \savesymbol{xLeftarrow}
   \savesymbol{xRightarrow} \savesymbol{xLeftrightarrow}
   \savesymbol{xrightharpoondown} \savesymbol{xrightharpoonup}
   \savesymbol{xleftharpoondown} \savesymbol{xleftharpoonup}
   \savesymbol{xleftrightharpoons} \savesymbol{xrightleftharpoons}
   \savesymbol{xhookleftarrow} \savesymbol{xhookrightarrow}
   \savesymbol{xmapsto} \savesymbol{underbracket}
   \savesymbol{overbracket} \savesymbol{lparen} \savesymbol{rparen}
   \savesymbol{dblcolon} \savesymbol{coloneqq} \savesymbol{Coloneqq}
   \savesymbol{coloneq} \savesymbol{Coloneq} \savesymbol{eqqcolon}
   \savesymbol{Eqqcolon} \savesymbol{eqcolon} \savesymbol{Eqcolon}
   \savesymbol{colonapprox} \savesymbol{Colonapprox}
   \savesymbol{colonsim} \savesymbol{Colonsim} \savesymbol{overbrace}
   \savesymbol{underbrace}

   \let\origAtBeginDocument=\AtBeginDocument
   \def\AtBeginDocument##1{##1}
   \usepackage[donotfixamsmathbugs]{mathtools}
   \let\AtBeginDocument=\origAtBeginDocument

   \restoresymbol{MTOOLS}{xleftrightarrow}
   \restoresymbol{MTOOLS}{xLeftarrow}
   \restoresymbol{MTOOLS}{xRightarrow}
   \restoresymbol{MTOOLS}{xLeftrightarrow}
   \restoresymbol{MTOOLS}{xrightharpoondown}
   \restoresymbol{MTOOLS}{xrightharpoonup}
   \restoresymbol{MTOOLS}{xleftharpoondown}
   \restoresymbol{MTOOLS}{xleftharpoonup}
   \restoresymbol{MTOOLS}{xleftrightharpoons}
   \restoresymbol{MTOOLS}{xrightleftharpoons}
   \restoresymbol{MTOOLS}{xhookleftarrow}
   \restoresymbol{MTOOLS}{xhookrightarrow}
   \restoresymbol{MTOOLS}{xmapsto}
   \restoresymbol{MTOOLS}{underbracket}
   \restoresymbol{MTOOLS}{overbracket} \restoresymbol{MTOOLS}{lparen}
   \restoresymbol{MTOOLS}{rparen} \restoresymbol{MTOOLS}{dblcolon}
   \restoresymbol{MTOOLS}{coloneqq} \restoresymbol{MTOOLS}{Coloneqq}
   \restoresymbol{MTOOLS}{coloneq} \restoresymbol{MTOOLS}{Coloneq}
   \restoresymbol{MTOOLS}{eqqcolon} \restoresymbol{MTOOLS}{Eqqcolon}
   \restoresymbol{MTOOLS}{eqcolon} \restoresymbol{MTOOLS}{Eqcolon}
   \restoresymbol{MTOOLS}{colonapprox}
   \restoresymbol{MTOOLS}{Colonapprox}
   \restoresymbol{MTOOLS}{colonsim} \restoresymbol{MTOOLS}{Colonsim}
   \restoresymbol{MTOOLS}{overbrace} \restoresymbol{MTOOLS}{underbrace}

  }
  {}
  \newif\ifXPFEIL
	
	\IfStyFileExists{extpfeil}
	  {\XPFEILtrue
	   \let\origRequirePackage=\RequirePackage
	   \renewcommand*{\RequirePackage}[2][]{}
	   \savesymbol{xlongequal}
	   \savesymbol{xmapsto}
	   \savesymbol{xhookrightarrow}
	   \usepackage{extpfeil}
	   \restoresymbol{XPFEIL}{xlongequal}
	   \restoresymbol{XPFEIL}{xmapsto}	   
	   \let\RequirePackage=\origRequirePackage
	  }
	  {}

	\usepackage{hyperref}
	\hyperbaseurl{.}	

\newtheorem*{theorem*}{Theorem}
\ifdefined\thmBookNumbering
	\newtheorem{theorem}{Theorem}[chapter]
	\newtheorem{definition}[theorem]{Definition}
	\newtheorem{remark}[theorem]{Remark}
	\newtheorem{corollary}[theorem]{Corollary}
	\newtheorem{proposition}[theorem]{Proposition}
	\newtheorem{lemma}[theorem]{Lemma}
	\newtheorem{example}[theorem]{Example}
	\newtheorem*{definition*}{Definition}
	\newtheorem*{maintheorem}{Main Theorem}
	\newtheorem*{unnumberedtheorem}{Theorem}
	\newtheorem*{remark*}{Remark}
	\newtheorem*{unnumberedcorollary}{Corollary}
\else
  \newtheorem{definition}{Definition}
	\newtheorem*{definition*}{Definition}
	\newtheorem{theorem}{Theorem}
	
	\newtheorem*{unnumberedtheorem}{Theorem}
	\newtheorem{remark}{Remark}
	\newtheorem*{remark*}{Remark}
	\newtheorem{corollary}{Corollary}
	
	\newtheorem{proposition}{Proposition}
	\newtheorem{lemma}{Lemma}
	\newtheorem{example}{Example}
\fi

	\newcommand{\grad}{\operatorname{grad}}
	\newcommand{\Div}{\operatorname{div}}
	\newcommand{\spann}{\operatorname{span}}

	\newcommand{\trace}{\operatorname{trace}}

	\newcommand{\pr}{\operatorname{pr}}	
	\newcommand{\To}{\longrightarrow}
	
	\newcommand{\rank}{\operatorname{rank}}	
	
	\newcommand{\im}{\operatorname{im}}
	\newcommand{\Hol}{\operatorname{Hol}}
	\newcommand{\hol}{\mathfrak{hol}}

	\newcommand{\LAkill}{\mathfrak{kill}}

	\newcommand{\coker}{\operatorname{coker}}
	\newcommand{\rmd}{{\rm d}}
	\newcommand{\rmi}{{\rm i}}
	\newcommand{\rme}{{\rm e}}
	\newcommand{\e}{\rme}

	\newcommand{\R}{\mathbb{R}}
	\newcommand{\PP}{\mathbb{P}}
	
	\newcommand{\Z}{\mathbb{Z}}
	
	\newcommand{\C}{\mathbb{C}}

	\newcommand{\T}{\mathbb{T}}	
		
	\newcommand{\IS}{\mathbb{S}}
	\newcommand{\IL}{\mathbb{L}}
	
	\newcommand{\mfX}{\mathfrak{X}}

	\newcommand{\M}{\mathcal{M}}
	
	\newcommand{\cN}{\mathcal{N}}
	\newcommand{\cM}{\mathcal{M}}	
	\newcommand{\cB}{\mathcal{B}}

	\newcommand{\cS}{\mathcal{S}}
	
	\newcommand{\tem}{\widetilde{\M}}

	\newcommand{\KNP}{\owedge}

	\newcommand{\Ric}{\operatorname{Ric}}

	\newcommand{\Hess}{\operatorname{Hess}}

\begin{document}

	\title{\textsc{On the Geometry of Circle Bundles with Special Holonomy}}
	  
  \author{\textsc{Daniel Schliebner}\thanks{The author is funded by the Berlin Mathematical School (BMS).}}
  \date{\today}

	\maketitle
	
  \begin{abstract}
  \noindent We investigate geometric properties of indecomposable but non-irreducible Lorentzian manifolds, 
  which are total spaces of circle bundles. We investigate under which conditions these manifolds are
  complete and give examples which fulfill the obtained conditions. In particular we investigate
  the Einstein equation for these spaces yielding  examples for complete compact Ricci flat Lorentzian manifolds 
  and manifolds with timelike Killing vector fields.
  Finally we study their holonomy and obtain in particular complete examples for
  Lorentzian manifolds with holonomy of so called type 4.
  
  \medskip
  
  
  \end{abstract}

	

	\section{Introduction}
\label{section:Intr}

In this paper we study certain Lorentzian manifolds\footnote{In this paper, all manifolds are assumed to be connected, smooth and without boundary.} $(\cM^{n + 2}, g)$ with \textit{special holonomy}, by which we mean that their holonomy representation $\rho : \Hol_x(\cM, g) \To {\rm O}(T_x\cM, g_x)$ for $x \in \cM$, acts indecomposable but non-irreducible, i.\,e.\ they admit a proper invariant degenerate subspace $W \subset T_x\cM$ but no proper non-degenerate subspace. Here, $\Hol_x(\cM,g) := \{ \mathcal P^g_\gamma \in {\rm O}(T_x\cM) \ | \ \gamma \text{ loop in } x \} \subset {\rm O}(T_x\cM)$ denotes the full \textit{holonomy group} of $(\cM, g)$ along piecewise smooth curves $\gamma$, closed in $x \in \cM$. Of course, since any connected Lie-subgroup $H \subset \operatorname{O}(1, n + 1)$ acting irreducibly on $\R^{1, n + 1}$ is equal to $\operatorname{SO}^0(1, n + 1)$ \cite{13}, any indecomposable Lorentzian manifold with restricted holonomy group $\Hol^0_x(\cM,g)$ (i.\,e.\ the subgroup of $\Hol_x(\cM,g)$ obtained by restricting to null-homotopic loops) not equal to $\operatorname{SO}^0(1, n + 1)$ has special holonomy. In particular, the Lorentzian manifolds with special holonomy play an important role within the classification of Lorentzian manifolds since, in the de Rham-Wu decomposition of any complete, simply-connected Lorentzian manifold, the Lorentzian factor is either $(\R, -dt^2)$, has holonomy $\operatorname{SO}^0(1, n + 1)$ or is a Lorentzian manifold with special holonomy. Beside this brief explanations, we refer, for example, to \cite{33} for a more comprehensive introduction. \par
	As in \cite[Section 2.2]{24} we will provide a sort of construction principle for Lorentzian manifolds $(\cM^{n + 2}, g)$ with recurrent or parallel light-like vector field $\xi \in \Gamma(T\cM)$ and hence special holonomy with non-trivial topology. The constructed manifolds $\cM$ are total spaces of $\IS^1$-bundles $\pi : \cM \To \cN$ with $c_1(\cM) = [\Psi]$ on which certain Lorentzian metrics $g = 2\rmi A \odot \pi^*\eta + f \cdot \pi^*\eta \odot \pi^*\eta + \pi^*h$ are defined, depending on a five-tuple of freely selectable objects $(A,\eta,f,\cN,h)$, where $h$ is a Riemannian metric on $\cN$ and $A \in \Omega^1(\cM, \rmi\R)$. We will refer to this construction by saying that \textit{$(\cM^{n + 2}, g)$ is of type $(\Psi, A, \eta, f)$ over $(\cN, h)$}. In \cite{24}, this construction was used to produce Lorentzian manifolds with special holonomy. In light of this it would be interesting if these constructions also yield special geometries such as (complete) Lorentzian Einstein spaces. As it turns out, the constructions of \textsc{L\"arz} with $\cN = \cB \times \IS^1$ and $\Psi \in \Omega^2(\cB)$ are complete (Theorem \ref{Prop:CompleteNonInt}) but unfortunately cannot be Einstein (Proposition \ref{Prop:RicciflatSphereBundle}). \par
	Another widely open field is to find examples for Lorentzian manifolds with special holonomy of so called \textit{type 3 or 4} (see Theorem \ref{prop:holTypes}) with special topological and geometric properties. 
	First examples were found in \cite{31,32,35} but to our knowledge it is unknown if there, for example, exist complete Lorentzian metrics of type 3 or 4. \par
	The main purpose of this paper is to study the Lorentzian manifolds $(\cM^{n + 2}, g)$ of type $(\Psi, A, \eta, f)$ over $(\cN, h)$ in light of the latter considerations. Namely, after Section \ref{section:Pre} with basic calculations to establish formulas for the Levi-Civita connection of $(\cM^{n + 2}, g)$, we provide different examples leading to complete (Section \ref{section:Indecomposables}) and Ricci-flat (Section \ref{Section:Geometry}) Lorentzian manifolds of type $(\Psi, A, \eta, f)$ over $(\cN, h)$:
	
\begin{unnumberedtheorem}[Theorem \ref{Thm:Ricciflat} + Corollary \ref{Cor:Ricciflat} + Corollary \ref{Cor:RicciflatComplete}]	
	Let $\cN := \cB \times \IS^1$ with $h := h_\cB \oplus du^2$ for an $n$-dimensional compact Ricci-flat Riemannian manifold $(\cB, h_\cB)$.
	Then for particular choices of $\Psi, A, \eta$ and $f$, the manifold $(\cM^{(n + 2)}, g)$ of type $(\Psi, A, \eta, f)$ over $(\cN, h)$ is
	a complete, compact and Ricci-flat Lorentzian manifold.
\end{unnumberedtheorem}
	
We stress that the constructions cannot provide examples for compact Lorentzian Einstein spaces with non-zero cosmological constant (Proposition \ref{Prop:RicciflatSphereBundle}). Moreover, these investigations also yield examples for compact Lorentzian manifolds with a global timelike Killing vector field, which must be complete, cf.\ \cite[Theorem 2.1]{15}. \par
	Section \ref{Section:Holonomy} is devoted to holonomy. By computing the universal cover of certain Lorentzian manifolds $(\cM^{n + 2}, g)$ of type $(\Psi, \eta, f)$ over $(\cN, h)$, cf.\ Proposition \ref{Thm:univCover}
and Proposition \ref{Thm:univCover2}, we can in particular compute the  holonomy of the complete, compact and Ricci-flat Lorentzian manifolds obtained above.

\begin{unnumberedtheorem}[Theorem \ref{Prop:HolUniv} + Corollary \ref{Prop:HolUniv2}]	
	Under a certain assumption on the fundamental group of $\cB$ (see Definition \ref{def:split}), the full holonomy group of the obtained complete, compact and Ricci-flat Lorentzian manifolds
	$(\cM^{(n + 2)}, g)$ of type $(\Psi, A, \eta, f)$ over $(\cB \times \IS^1, h)$ equals
	$\Hol(\cM^{(n + 2)}, g) = \Hol(\cB,h_\cB) \ltimes \R^n$.
\end{unnumberedtheorem}

	At the end of Section \ref{Section:Holonomy} we also give a construction principle for complete Lorentzian manifolds with holonomy of type 4 using
the Lorentzian manifolds $(\cM^{n + 2}, g)$ of type $(\Psi, A, \eta, f)$ over $(\cN, h)$ which are, to our best knowledge, the first examples of this kind.

\begin{unnumberedtheorem}[Theorem \ref{Thm:Type4CplExamples}]	
	For each Abelian Lie subalgebra $\mathfrak{g} \subset \mathfrak{so}(k)$ there exists a complete indecomposable Lorentzian manifold with holonomy of type 4 possessing $\mathfrak{g}$ as orthogonal part.
\end{unnumberedtheorem}
	
	\section{Preliminaries}
\label{section:Pre}

\subsection{Total Spaces of Circle Bundles}
\label{Subsection:totalSphere}

This section is devoted to the presentation of the construction of the stated Lorentzian metrics on the total spaces of $\IS^1$-bundles. \par
	Let $(\cN^{n + 1}, h)$ be an $(n + 1)$-dimensional Riemannian manifold and $\omega \in H^2(\cN, \Z)$. For the $\IS^1$-bundle $\pi : \cM \To \cN$ with first Chern class $c_1(\cM) = \omega$ consider the following Lorentzian metric $g$ on $\cM$. Take any closed 2-form $\overline{\Psi} \in \Omega^2(\cN)$ s.\,t.\ $\overline{\Psi}$ represents $\omega$ in the de Rham cohomology and a corresponding connection $A \in \Omega^1(\cM, \rmi\R)$ with curvature $F^A = dA = -2 \pi \rmi \pi^*\overline{\Psi}$. Then, for any nowhere vanishing closed 1-form $\eta \in \Omega^1(\cN)$ and
any function $f \in C^\infty(\cM)$ define
\begin{equation}
	\label{equ:Metric}
	g := 2 \rmi A \odot \pi^*\eta + f \cdot \pi^*\eta \odot \pi^*\eta + \pi^*h.
\end{equation}
Then, $(\cM^{n + 2}, g)$ is an $(n + 2)$-dimensional Lorentzian manifold. \par
	Henceforth, we write $\Psi := \pi \cdot \overline{\Psi}$, and thus $F^A = -2 \rmi \pi^*\Psi$. To refer to this construction we make the following definition.

\begin{definition}
	\label{def:SphereMF}
	The Lorentzian manifold $(\cM^{n + 2}, g)$ with $g$ chosen as in (\ref{equ:Metric}) is called \dStr{manifold of type $(\Psi, A, \eta, f)$ over $(\cN, h)$}.
\end{definition}

For the upcoming calculations we will use the following local frame on $(\cM, g)$. Let $x = \pi(y) \in \cN$ be an arbitrary point on $\cN$. On $\cN$ we have the global vector field $E_\eta := \frac{\eta^\sharp}{||\eta^\sharp||_h^2}$ and on $\cM$ the fundamental vector field $\xi \in \Gamma(T\cM)$ corresponding to the $\IS^1$-action, i.\,e.\ 
$$
	\xi(z) := \widetilde{\rmi}(z) = \frac{d}{dt}(z \cdot \exp(t \cdot \rmi))|_{t = 0},
$$
$z \in \cM$, which is light-like w.\,r.\,t.\ $g$. Locally around $x \in U \subset \cN$, we may choose a frame $E_1, \ldots, E_n, E_\eta$ s.\,t.\ $h(E_i,E_j) = \delta_{ij}$ and
$\ker \eta = \spann\{E_1, \ldots, E_n\} {\bot_h} \R E_\eta$. Taking its horizontal lifts $E_i^* \in \Gamma(T\cM|_{\pi^{-1}(U)})$ we thus obtain a local orthonormal frame on $(\cM, g)$:
\begin{equation}
	\label{equ:Basis}
	\rme_i := E_i^*,
	\ \ \ \ \rme_+ := \zeta + \tfrac{1}{2} H \xi,
	\ \ \ \ \rme_- := \rme_+ + \xi,
\end{equation}
with $\zeta := E_\eta^*$, $H := (f + \tfrac{1}{||\eta^\sharp||_h^2} - 1)$ and $i = 1, \ldots, n$. 
Then, $g(\rme_i, \rme_j) = \delta_{ij}$, $g(\rme_i, \rme_+) = g(\e_i,\rme_-) = 0$, $g(\rme_+,\rme_+) = 1$ and
$g(\rme_-, \rme_-) = -1$. \par
	We do now proceed to calculate the Levi-Civita connection corresponding to $g$. Note that in all forthcoming formulas, the Latin indices $i,j,k$ and $\ell$ run from $1$ to $n$ and $\xi$, $+$ denoted \textit{as index within tensors} means plugging in the vector field $\xi$ or $\rme_+$, respectively. Moreover, we omit the components with at least one $\rme_-$-vector since these are immediate by the multi-linearity and Leibniz-rules of the objects in question.

\begin{lemma}
	\label{equ:CovDer}
	Let $(\cM^{n + 2}, g)$ be of type $(\Psi, A, \eta, f)$ over $(\cN, h)$. Then,
	\begin{enumerate}
		\item[(a)] $\nabla^g_{\rme_i} \rme_j 		   = \overline{\nabla^h_{E_i} E_j}^* + (\tfrac{\rmi}{2}F^A(\rme_i, \rme_j) - h(E_\eta, \nabla^h_{E_i}E_j)) \xi$,
		\item[(b)] $\nabla^g_{\rme_+} \rme_j			 = \overline{\nabla^h_{E_\eta} E_j}^* + \overline{\psi(E_j)}^* - (\rmi F^A(\e_j, \e_+) + \tfrac{1}{2}dH(\e_j))\xi$,
		\item[(c)] $\nabla^g_{\rme_i} \rme_+ 		 	 = \overline{\nabla^h_{E_i} E_\eta}^* + \overline{\psi(E_i)}^*$,
		\item[(d)] $\nabla^g_{\rme_+} \rme_+ 		   = \overline{\nabla^h_{E_\eta} E_\eta}^* + 2\overline{\psi(E_\eta)}^* - \tfrac{1}{2}\grad_g \! f - \tfrac{1}{2}\e_+(f)\xi$,
		\item[(e)] $\nabla^g \xi 									 = -\tfrac{1}{2} \xi(f) \cdot \pi^*\eta \otimes \xi$.
	\end{enumerate}
	Here, $\psi \in \Omega^1(\cN, T\cN)$ is defined as $h(\psi(E_i), E_j) := \Psi(E_i, E_j)$ and for any $X \in \Gamma(T\cN)$ we define 
	$\overline{X} := \pr_{\ker \eta} X$. Hence, $X = \overline{X} + \eta(X)E_\eta$.
\end{lemma}

\begin{proof}
	Since $\e_i = E_i^*$ are horizontal lifts, one has $[\e_i, \xi] = [\zeta, \xi] = 0$ and 
	\begin{eqnarray*}
		\,[X^*, Y^*]   & = & [X, Y]^* - \widetilde{F^A(X, Y)} = [X,Y]^*  + \rmi F^A(X^*,Y^*)\xi, \\
		\,[\e_+, X^*]  & = & [\zeta, X^*] + \tfrac{1}{2} [H\xi, \e_i] = [E_\eta, X]^* + \rmi F^A(\zeta, X^*) - \tfrac{1}{2}dH(X^*) \xi,
	\end{eqnarray*}
	for all $X,Y \in \Gamma(T\cN)$. Moreover, by taking into account that $\eta \in \Omega^1(\cN)$ is closed and $\eta(E_\eta) \equiv 1$, we see that
	$$
		\eta([E_i, E_j]) = \eta([E_\eta, E_i]) = 0
	$$
	for $i = 1, \ldots, n$.	The formulas \textit{(a)} to \textit{(e)} are now immediate consequences of the Koszul formula for $\nabla^g$.
\end{proof}

\subsection{The Canonical Screen bundle}
\label{Subsection:screen}

Having special holonomy, the Lorentzian manifolds $(\cM^{n + 2}, g)$ of type $(\Psi, A, \eta, f)$ over $(\cN, h)$ admit a holonomy invariant null line $L := W \cap W^\bot$ giving rise to a line bundle $\IL \subset T\cM$. Being aware of the inclusion $\IL \subset \IL^\bot$, one naturally gains a flag
\begin{equation}
\label{equ:flag}
\IL \subset \IL^\bot \subset T\cM
\end{equation}
and an $n$-dimensional quotient bundle $\Sigma := \operatorname{coker}(\IL \hookrightarrow \IL^\bot) = \IL^\bot / \IL$ (the \textit{screen bundle}). Taking into account a splitting $s : \Sigma \To \IL^\bot$ of the exact sequence $0 \to \IL \hookrightarrow \IL^\bot \twoheadrightarrow \Sigma \to 0$, one obtains an $n$-dimensional distribution
$\IS := \im s$ which is called a \textit{screen distribution} of $(\cM^{n + 2}, g)$. Note that, as we will see, possible consequences on the geometry or topology of $(\cM^{n + 2}, g)$ depend on the existence of certain realizations of $\IS$. A realization $\IS$ is called \textit{horizontal} if $[\Gamma(\IL), \Gamma(\IS)] \subset \Gamma(\IS)$ and \textit{integrable} if so is the distribution $\IS$ itself. If one finds horizontal or integrable realizations of $\IS$ this turns out to be very useful, for example these properties were used in \cite{21,24,schHigh1stbetti} to prove several results concerning topology and geometry. \par
	In the case of the circle bundle metrics studied in this paper we clearly have that $\IL = \R\xi$ and, locally, $\IL^\bot = \spann \{\xi, \e_1, \ldots, \e_n\}$. Moreover we have a canonical realization of the screen bundle. Namely, we may define by
\begin{equation}
	\label{equ:Z}
	Z := \frac{1}{2}\xi - \e_-
\end{equation}
a light-like vector field with $g(\xi, Z) = 1$. Then, the metric $g$ is non-degenerate on the plane $\spann\{ \xi, Z \}$ and we obtain a realization 
of the screen bundle by $\IS := \spann\{\xi, Z\}^{\bot_g}$ with nice properties:

\begin{lemma}
	\label{lemma:Screen}
	Let $(\cM^{n + 2}, g)$ be of type $(\Psi, A, \eta, f)$ over $(\cN, h)$. Then, realizing the screen bundle as $\IS = \{\xi, Z\}^{\bot_g}$,
	we obtain a horizontal realization of the screen bundle. Moreover, the screen distribution $\IS$ is integrable if and only if 
	$F^A|_{\ker \pi^*\eta \times \ker \pi^*\eta} = 0$ or, equivalently, $\eta \wedge \Psi = 0$.
\end{lemma}

\begin{proof}
	Of course, choosing, locally, the orthonormal frame (\ref{equ:Basis}), we clearly have that $\IS|_{\pi^{-1}(U)} = \spann\{\e_1, \ldots, \e_n\}$. Now, since 
	$[\xi, \e_i] = 0$, $\IS$ is horizontal. Moreover, we have 
	$$
		[\e_i, \e_j] = [E_i, E_j]^* + \rmi F^A(\e_i, \e_j)\xi.
	$$
	Consequently, $[\e_i, \e_j] \in \Gamma(\IS)$ if and only if $F^A(\e_i, \e_j) = 0$ or, equivalently, $\eta \wedge \Psi = 0$.
\end{proof}
	
	Indeed, this does in general not imply that we cannot find another realization of the screen bundle which is integrable and horizontal.
But in fact, for the following manifold of type $(\Psi, A, \eta, f)$ over $(\cN, h)$ one can prove that such a realization cannot exist.

\begin{example}[{\cite[Example 1]{21}}]
	\label{Ex:torusNonIH}
	Let $\cN = \T^n \times \IS^1$ with the metric $h = h_{\text{flat}} \oplus du^2$ and $0 \neq \omega \in H^2(\T^n, \Z) \cap H_{\rm dR}^2(\T^n)$. 
	Then for $\eta := du$, $\Psi \in \omega$ and any $f \in C^\infty(\cM)$, for the manifold $(\cM^{n + 2}, g)$ of type $(\Psi, A, \eta, f)$ over $(\T^n \times \IS^1, h)$ there exists no integrable realization
	of the screen bundle.
\end{example}

	\section{Completeness}
\label{section:Indecomposables}

We are now interested in conditions for which the Lorentzian manifolds $(\cM^{n + 2}, g)$ of type $(\Psi, A, \eta, f)$ over $(\cN, h)$ are complete. To establish criteria for completeness we preliminarily prove the following proposition which is a slight generalization of \cite[Proposition 2.1]{22} in the Lorentzian case. However, for the sake of completeness, we present the proof here.

\begin{proposition}
	\label{CompleteSanch}
	Let $(\cM^{n + 2}, g)$ be a Lorentzian manifold with timelike vector field $X$ that satisfies the following three conditions:	
	\begin{itemize}
		\item[(i)]   $g(X,X)^{-1}$ is bounded on $\cM$,		
		\item[(ii)] the Riemannian metric $g^R$ given by
				$$
					g^R|_{X^\bot \times X^\bot} = g, \ g^R(X,X) = -g(X,X), \ g^R|_{X^\bot \times X} = g^R|_{X \times X^\bot} = 0,
				$$
				is complete.		
	\end{itemize}
	Then for every inextensible $g$-geodesic  $\gamma : [0, \varepsilon) \To \cM$, the map 
	$$
		t \in [0, \varepsilon) \longmapsto (\mathscr L_X g)(\dot\gamma(t), \dot\gamma(t))
	$$
	is unbounded.	(Here $\mathscr L_X g$ denotes the Lie-derivative of $g$ along $X$.)
\end{proposition}

\begin{proof}
	Let $\gamma : [0, \varepsilon) \To \cM$ be an inextensible $g$-geodesic with $0 < \varepsilon < \infty$. It suffices to show that the function
	$t \in [0, \varepsilon) \longmapsto g^R(\dot\gamma(t), \dot\gamma(t)) \in \R$ is bounded. Namely, in this case, $\{x_n := \gamma(t_n)\}$ for some $\{t_n\} \rightarrow \varepsilon$ 
	is a $d_R$-Cauchy sequence, where $d_R$ denotes the geodesic distance w.r.t.\ $g^R$. Since $g^R$ is complete, the closure of $\{x_n\}$ 
	is compact and so there exists a convergent subsequence to, say, $x \in M$. But as $\{x_n\}$ is Cauchy, it converges to $x$, too, while the sequence
	$\{t_n\}$ with $t_n \rightarrow \varepsilon$ can be chosen arbitrarily. But then $\gamma : [0, \varepsilon) \To \cM$ is extensible beyond $\varepsilon$ via
	$\lim_{t \rightarrow \varepsilon^-}\gamma(t) := x$ which is a contradiction. \par
		Let $\widehat{X} := X/||X||$. Since $g(\widehat{X}, \widehat{X}) = -1$, $g^R(\widehat{X}, \widehat{X}) = 1$ 
	and $\pr_{\widehat{X}^\bot} \dot\gamma = \dot\gamma + g(\widehat{X},\dot\gamma)\widehat{X}$, we obtain
	\begin{eqnarray*}
		g(\dot\gamma, \dot\gamma) 
			& = & g(\pr_{\widehat{X}^\bot} \dot\gamma - g(\widehat{X},\dot\gamma)\widehat{X}, \pr_{\widehat{X}^\bot} \dot\gamma - g(\widehat{X},\dot\gamma)\widehat{X}) \\
			& = & g(\pr_{\widehat{X}^\bot} \dot\gamma, \pr_{\widehat{X}^\bot} \dot\gamma) - 2g(\pr_{\widehat{X}^\bot} \dot\gamma, g(\widehat{X},\dot\gamma)\widehat{X})
						 + g(\widehat{X},\dot\gamma)^2g(\widehat{X}, \widehat{X}) \\
			& = & g(\pr_{\widehat{X}^\bot} \dot\gamma, \pr_{\widehat{X}^\bot} \dot\gamma) - g(\widehat{X},\dot\gamma)^2 \\
			\text{ and } \\
	  g^R(\dot\gamma, \dot\gamma) 
	  	& = & g^R(\pr_{\widehat{X}^\bot} \dot\gamma - g(\widehat{X},\dot\gamma)\widehat{X}, \pr_{\widehat{X}^\bot} \dot\gamma - g(\widehat{X},\dot\gamma)\widehat{X}) \\
	  	& = & g^R(\pr_{\widehat{X}^\bot} \dot\gamma, \pr_{\widehat{X}^\bot} \dot\gamma) - 2g^R(\pr_{\widehat{X}^\bot} \dot\gamma, g(\widehat{X},\dot\gamma)\widehat{X})
						 + g(\widehat{X},\dot\gamma)^2g^R(\widehat{X}, \widehat{X}) \\
			& = & g^R(\pr_{\widehat{X}^\bot} \dot\gamma, \pr_{\widehat{X}^\bot} \dot\gamma) + g(\widehat{X},\dot\gamma)^2.
	\end{eqnarray*}
	Since $g^R|_{X^\bot \times X^\bot} = g$ it follows
	$$
		g^R(\dot\gamma, \dot\gamma) = g(\dot\gamma, \dot\gamma) + \frac{2}{g(X,X)} g(X,\dot\gamma)^2.
	$$
	Since $g(\dot\gamma, \dot\gamma)$ is constant and $g(X,X)^{-1}$ is bounded, we are left to show that $g(X,\dot\gamma)$ is bounded on $[0, \varepsilon)$. \par
		We compute
	$$
		\frac{d}{dt}g(X,\dot\gamma) = \frac{1}{2}(\mathscr L_X g)(\dot\gamma,\dot\gamma).
	$$
	Hence, if $(\mathscr L_X g)(\dot\gamma,\dot\gamma)$ is bounded, so is $\frac{d}{dt}g(X,\dot\gamma)$ and consequently, also $g(X,\dot\gamma)$ on $[0, \varepsilon)$.
\end{proof}

\noindent With the aid of the former proposition we can now prove the following.

\begin{theorem}
	\label{Prop:CompleteNonInt}
	Let $(\cM^{n + 2}, g)$ be of type $(\Psi, A, \eta, f)$ over $(\cN, h)$ with compact base $\cN$ s.\,t.\ the function $f \in C^\infty(\cM)$ 
	is constant along the fibers and $\eta^\sharp$ is a Killing field on $(\cN,h)$. 
	If, moreover, $\Psi(\eta^\sharp, \cdot) = 0$ then $(\cM^{n + 2}, g)$ is complete. \par
		In particular, on $\cM$ there exists a nowhere-vanishing timelike Killing vector field if it additionally holds $\zeta(f) = 0$.
\end{theorem}

\begin{proof}
	We define a vector field $K \in \Gamma(T\cM)$ by
	\begin{equation}
		\label{equ:killing}
		K := \zeta + \tfrac{C}{2} \cdot \xi \ \ \text{ with constant } \ \ C := \max\limits_{\cM} g(\zeta, \zeta) + \varepsilon \in \R,
	\end{equation}
	where $\varepsilon > 0$ is arbitrary chosen. \par
		To apply Proposition \ref{CompleteSanch}, we have to show that $K \in \Gamma(T\cM)$ is timelike as the conditions \textit{(i)} and \textit{(ii)} are satisfied since $\cM$ is compact. \par				
	For the length of $K$ we get
	$$
		g(K, K) < -\varepsilon < 0
	$$
	due to the definition of $C \in \R$. \par
		Since $\Psi(\eta^\sharp, \cdot) = 0$, we obtain by the formulas in Lemma \ref{equ:CovDer}, and the fact that
	$$
		\nabla^g K = \nabla^g \zeta = \nabla^g \e_+ - \tfrac{1}{2}dH \otimes \xi
	$$	
	(since $\nabla^g \xi = 0$) that the Lie-derivative $\mathscr L_K g$ is given by
	$$
		\mathscr L_K g = \frac{1}{2} \zeta(f) \pi^*\eta \odot \pi^*\eta.
	$$
	Since $\pi^*\eta = -g(\xi, \cdot)$ is $\nabla^g$-parallel and $\zeta(f)$ is bounded as $\cM$ is compact, there is no inextensible geodesic on $(\cM, g)$ by Proposition \ref{CompleteSanch}, hence
	completeness follows.
\end{proof}

As we will see in the next section there are quite a lot of examples that fulfill the assumptions made in the previous theorem and are hence geodesically complete. Of course, the assumption
$\Psi(\eta^\sharp, \cdot) = 0$ is not absolutely necessary. Indeed, the next proposition gives examples for compact manifolds of type $(\Psi, A, \eta, f)$ over $(\cN, h)$ with compact base $\cN$
and $\Psi(\eta^\sharp, \cdot) \neq 0$ which are complete, too.

\begin{proposition}
	\label{Prop:CompleteInt}
	Let $(\cM^{n + 2}, g)$ be of type $(\Psi, A, \eta, f)$ over $(\cN, h)$ with compact base $\cN$ s.\,t.\ the function $f \in C^\infty(\cM)$ 
	is constant along the fibers. Let either	
	\begin{itemize}
		\item[(i)] $\alpha \in \Omega^1(\cN)$ be a $h$-parallel 1-form with $\eta(\alpha^\sharp) = 0$ or
		\item[(ii)] $\cN = \cB \times \IS^1$, $\eta = du$ the coordinate 1-form on $\IS^1$, $b_1(\cB) = 0$ and $\alpha$ a closed 1-form on $\cB$.
	\end{itemize}
	Then, choosing $\Psi := \alpha \wedge \eta$, the manifold $(\cM^{n + 2}, g)$ is complete.
\end{proposition}

\begin{proof}
	Let $K \in \Gamma(T\cM)$ as in the proof before. In this case, we have that
	$$
		\mathscr L_K g = \frac{1}{2} \zeta(f) \pi^*\eta \odot \pi^*\eta + 2\pi^*\alpha \odot \pi^*\eta.
	$$
	Assume there is an inextensible geodesic $\gamma : [0,\varepsilon) \To \cM$. To prove \textit{(ii)}, let $\alpha = d\tau$. If we denote by $\delta := \pr_\cB \circ \pi \circ \gamma$
	the projected curve on $\cB$, then by Lemma \ref{equ:CovDer}, $\delta$ is a $h_\cB$-geodesic and as $\tau \in C^\infty(\cB)$,
	$$
		\pi^*\alpha(\dot\gamma) = d\tau(d\pi(\dot\gamma)) = d\tau(\dot\delta) = h_\cB(\grad_{h_\cB} \tau, \dot\delta) \leq ||\grad_{h_\cB} \tau|| \cdot ||\dot\delta||
	$$
	by the Cauchy-Schwarz inequality. Hence, $\pi^*\alpha(\dot\gamma)$ is bounded. Since $\cN$ and $\cM$ are compact, $\zeta(f)$ is bounded, while $\pi^*\eta = -g(\xi, \cdot)$ is $\nabla^g$-parallel.
	Hence, $(\mathscr L_K g)(\dot\gamma, \dot\gamma)$ is bounded and the assertion now follows from Proposition \ref{CompleteSanch}. \par	
		For case \textit{(i)}, by $\nabla^h\alpha^\sharp = 0$ and the formula in Lemma \ref{equ:CovDer}b), we see that
	$$
		\frac{d}{dt}\pi^*\alpha(\dot\gamma) 
		= \frac{d}{dt}g(\pi^*\alpha^\sharp, \dot\gamma(t))
		= (2\alpha(\alpha^\sharp)\eta(\eta^\sharp) + \tfrac{1}{2}df(\alpha^\sharp) + 1) \cdot \pi^*\eta(\dot\gamma).		
	$$
	Since $\eta(\eta^\sharp), \alpha(\alpha^\sharp)$ and $df(\alpha^\sharp)$ are bounded $\frac{d}{dt}\alpha(\dot\gamma)$ is bounded on $[0, \varepsilon)$ and hence so is again $\pi^*\alpha(\dot\gamma)$.
	The same arguments as in \textit{(i)} complete the proof.
\end{proof}

As we have already mentioned in the introduction, the Lorentzian manifolds of type $(\Psi, A, \eta, f)$ over $(\cN, h)$ were already 
studied in \cite{24} to produce Lorentzian manifolds with special holonomy. Namely, there it was proven that for particular choices of $(\Psi, \eta, f)$ and the base manifold $(\cN, h)$, the resulting manifolds $(\cM^{n + 2}, g)$ of type $(\Psi, A, \eta, f)$ over $(\cN, h)$ have full holonomy
$\Hol(\cM^{n + 2}, g) = (\R^+ \times G) \ltimes \R^n$ resp. $\Hol(\cM^{n + 2}, g) = G \ltimes \R^n$ for recurrent resp. parallel fundamental vector field $\xi \in \Gamma(T\cM)$, where $G := \Hol(\cB, h_\cB)$ is the holonomy group of a certain Riemannian manifold $(\cB, h_\cB)$ and $\cN = \cB \times \IS^1$.\footnote{In \cite{24}, the manifold $(\cM^{n + 2}, g)$ then is called of \textit{toric type}.} 
In particular, in \cite[Prop. 2.42]{24} it is proven that taking $\cN = \T^n = \T^{n - 1} \times \IS^1$ and $\Psi = du \wedge dv$, the resulting manifold $(\cM^{n + 2}, g)$ of type $(du \wedge dv, A, du, f)$ over $(\T^n, g_{\T^n})$ is complete. This result however turns out to be a special case of our Proposition \ref{Prop:CompleteInt}. Moreover, all provided compact examples with special holonomy and base $\cN = \cB \times \IS^1$ in \cite{24} are complete by Theorem \ref{Prop:CompleteNonInt}, when $f \in C^\infty(\cM)$ is chosen to be constant along the fibers.
	
	\section{Geometry}
\label{Section:Geometry}

A possible question in the discussed construction is, whether the obtained Lorentzian manifolds with special holonomy produce examples with certain distinguished geometries. Indeed, for the case $\cM = \R^n$, 
a similar family of metrics was studied in \cite{17,18}, where \textsc{Gibbons}, \textsc{Pope}, \textsc{Leistner} and \textsc{Galaev} considered conditions under which certain Walker metrics produce Einstein metrics. An in some sense generalized but global version of the Walker metrics they considered is the presented construction of Lorentzian manifolds $(\cM^{n + 2}, g)$ of type $(\Psi, A, \eta, f)$ over $(N, h)$. The present section therefore deals with the question, whether these constructions produce Ricci-flat or even Einstein metrics with non-zero cosmological constant. As it turns out, the former is possible, while the latter is not due to the fact that the Hessian of $f \in C^\infty(\cM)$ cannot be constant on $\xi \times \xi$.
Together with the former considerations in this paper we thus additionally obtain completeness results for the obtained Ricci-flat Lorentzian manifolds. 

We proceed to present the formulas for the Riemannian curvature tensor $\mathcal R^g$ and the Ricci tensor $\Ric^g$, where we use the sign convention
$$
	\mathcal R^g(X,Y)Z := \nabla^g_X\nabla^g_Y Z -  \nabla^g_Y\nabla^g_X Z -  \nabla^g_{[X,Y]} Z.
$$
With the symbol $\KNP$ we denote the Kulkarni–Nomizu product.

\begin{lemma}
	\label{lemma:Curvature}
	Let $(\cM^{n + 2}, g)$ be of type	$(\Psi, A, \eta, f)$ over $(\cN, h)$. Then the only non-vanishing terms of $\mathcal R^g$ are the following:
	\begin{eqnarray}
		\mathcal R^g_{ijk\ell} 
			& = & \mathcal R^h_{ijk\ell} + \tfrac{1}{||\eta^\sharp||^2} (\nabla^h \eta \KNP \nabla^h \eta)(E_i,E_j,E_k,E_\ell), \\
		\mathcal R^g_{i++j} 
			& = & \mathcal R^h_{i \eta \eta j} + \tfrac{1}{||\eta^\sharp||^2} (\nabla^h \eta \KNP \nabla^h \eta)(E_i,E_\eta,E_\eta,E_j) \nonumber \\
			&   & + \ 2(\Psi(\cdot, E_\eta) \odot (\nabla^h \eta)(E_\eta))(E_i, E_j) + (\nabla^h_{E_i} \Psi)(E_\eta, E_j) + (\nabla^h_{E_j} \Psi)(E_\eta, E_i) \nonumber \\
			&   & + \ h(\overline{\psi(E_i)}, \overline{\psi(E_j)}) - \tfrac{1}{2}(\Hess_gf)(\e_i, \e_j), \\
		\mathcal R^g_{ijk+} 
			& = & \mathcal R^h_{ijk\eta} + \tfrac{1}{||\eta^\sharp||^2} (\nabla^h \eta \KNP \nabla^h \eta)(E_i,E_j,E_k,E_\eta)\nonumber \\
			&   & + \ (\Psi(\cdot, E_\eta) \wedge (\nabla^h \eta)(E_k))(E_i, E_j) + (\nabla^h_{E_k} \Psi)(E_i, E_j), \\
		\mathcal R^g_{i++\xi} 
			& = &  -\tfrac{1}{2} (\Hess_g f) (\e_i, \xi), \\
		\mathcal R^g_{+\xi\xi+} 
			& = & -\tfrac{1}{2} (\Hess_g f) (\xi, \xi).
	\end{eqnarray}	
\end{lemma}

\begin{proof}
	The proof is straightforward by Lemma \ref{equ:CovDer}. Note that the $(2,0)$-tensor $\nabla \eta$ is symmetric since $\eta$ is closed. 
	Namely, as $0 = d\eta(X,Y) = X(\eta(Y)) - Y(\eta(X)) - \eta([X,Y])$, one infers
	$$
		(\nabla^h_X \eta)(Y) = X(\eta(Y)) - \eta(\nabla^h_XY) \stackrel{d\eta = 0}{=} Y(\eta(X)) + \eta([X,Y]) - \eta(\nabla^h_XY) = (\nabla^h_Y \eta)(X),
	$$
	which justifies the term $(\nabla^h \eta \KNP \nabla^h \eta)$. Moreover, $\Psi$ satisfies the second Bianchi identity since it is closed.
\end{proof}

To make notation short, we define the symmetric tensor $T_\eta : \Gamma(T\cM) \times \Gamma(T\cM) \To \R$ as follows:
\begin{equation}
	\label{equ:Teta}
	T_\eta(X,Y) := \frac{1}{||\eta^\sharp||^2}\sum_{k = 1}^n(\nabla^h \eta \KNP \nabla^h \eta)(d\pi(X), \e_k, \e_k, d\pi(Y)).
\end{equation}
By contraction of $\mathcal R^g$ we infer the non-vanishing terms of the Ricci tensor.
\begin{lemma}
	\label{lemma:Ric}
	Let $(\cM^{n + 2}, g)$ be of type	$(\Psi, A, \eta, f)$ over $(\cN, h)$. Then the only non-vanishing terms of $\Ric^g$ are the following:
	\begin{eqnarray}
		\Ric^g_{ij}
			& = & \Ric^h_{ij} + T_\eta(\e_i,\e_j), \\
		\Ric^g_{i+}
			& = & \Ric^h_{i\eta} + T_\eta(\e_i,\e_+) + \tfrac{\rmi}{2}((\Div_g F^A)(\e_i) + F^A(\e_i, \zeta)\Div_g\pi^*\eta) \nonumber \\
			&   & - \ \Psi(\nabla^h \eta^\sharp, E_\eta) - \tfrac{1}{2}(\Hess_g f)(\e_i, \xi), \\
		\Ric^g_{++}
			& = & \Ric^h_{\eta\eta} + T_\eta(\e_+,\e_+) + 2\trace_h[\Psi(\cdot, E_\eta) \odot (\nabla^h \eta)(E_\eta)] \nonumber \\
			&   & - \ 2 (\Div_h \Psi)(E_\eta) + ||\overline{\psi}||_h - \tfrac{1}{2}\Delta_g f - \xi(f)(\tfrac{1}{2}\xi(f) + 1),
			\label{ref:equRic2} \\
		\Ric^g_{\xi+}
			& = & - \tfrac{1}{2} (\Hess_g f) (\xi, \xi),
			\label{ref:equRic1}
	\end{eqnarray}
	where $\Div$ is the divergence of a tensor.\footnote{Let $T$ be a $(r,0)$ tensor and $g_0$ a semi-Riemannian metric. Then we define the divergence of $T$ by 
	$\Div_{g_0}T := \sum_{k} \varepsilon_k (\nabla_{e_k} T)(e_k, \cdot, \ldots, \cdot)$, where $e_i$ is a $g_0$-orthonormal frame with $\varepsilon_k := g_0(e_k,e_k)$.}
\end{lemma}

\begin{remark}
	\label{rem:Vanish}
	If $\eta \in \Omega^1(\cN)$ is recurrent, i.\,e.\ $\nabla^h \eta = \alpha \otimes \eta$ for some $\alpha \in \Omega^1(\cN)$, then
	$(\nabla^h \eta \KNP \nabla^h \eta)$ and hence $T_\eta$ already vanishes identically.
\end{remark}

\begin{remark}
	\label{rem:GibPope}
	To compare the curvature equations of Lemma \ref{lemma:Curvature} and Lemma \ref{lemma:Ric} with the results in \cite{17}, note that in their notation,
	$F_{\alpha\beta} = -F_{\alpha\beta}^A$ and $g_{\alpha\beta} = \delta_{\alpha\beta} \equiv const$.
\end{remark}

\noindent As the following theorem proves, this construction yields examples for Ricci-flat manifolds, even in the non-trivial case where $(\cN, h)$ is Ricci-flat but $\Psi \neq 0$. An obvious obstruction is the fact that for $g$ to be Ricci-flat, $f \in C^\infty(\cM)$ must be constant along the fibers due to (\ref{ref:equRic1}).

\begin{theorem}
	\label{Thm:Ricciflat}
	Let $\cN := \cB \times \IS^1$ or $\cN := \cB \times \R$ with $h := h_\cB \oplus du^2$ for an $n$-dimensional Riemannian manifold $(\cB, h_\cB)$.
	Moreover, let $(\cB, h_\cB)$ be Ricci-flat and $\eta := du$. Choose $\omega \in H_{\rm dR}^1(\cB) \cap H^1(\cB,\Z)$ and a representative 
	$\alpha \in \omega$ and consider the $\IS^1$-bundle $\pi : \cM \To \cN$ with $c_1(\cM) = [\alpha \wedge \eta]$. 
	Finally, choose $\Psi := \alpha \wedge \eta$ and $f := \widehat{f} \circ \pi \in C^\infty(\cM)$, where $\widehat{f} := f_\cB \cdot f_{\IS^1}$ with $f_\cB \in C^\infty(\cB)$ and 
	$f_{\IS^1} \in C^\infty(\IS^1)$. \par	
		Then, the Lorentzian manifold $(\cM^{(n + 2)}, g)$ of type $(\Psi, A, \eta, f)$ over $(\cN, h)$ is Ricci-flat if and only if
	$\Delta_{h_\cB}(f_\cB) = -4\Div_{h_\cB}(\alpha)$.
\end{theorem}

\begin{proof}
	Due to the definition of $h$ and $\eta$, $\nabla^h \eta = 0$. As $\alpha$ and $\eta$ are linearly independent, 
	we may choose on $\cB$ a local orthonormal frame $E_1, \ldots, E_n \in \mfX(U), \ U \subset \cB$ and consider the corresponding basis
	as in (\ref{equ:Basis}). Therefore, $F^A(\e_i,\e_j) = 0$ and $\Psi(E_i, E_j) = 0$ for all $i,j = 1, \ldots, n$. We obtain:
	$$
		(\Div_h \Psi)(E_\eta) = \Div_{h_\cB}(\alpha) = \Div_{h_\cB}(\alpha^\sharp).
	$$
	As $(\cB, h_\cB)$ is Ricci-flat, (\ref{ref:equRic2}) turns into 
	$$
		\Ric^g_{++} = -\frac{1}{2}\Delta_{h_\cB}(f_\cB) - 2\Div_{h_\cB}(\alpha),
	$$
	which proves the theorem.
\end{proof}

\noindent For the existence of concrete examples one needs to find solutions of the Poisson equation 
$$
	\Delta_{h_\cB}(f_\cB) = -4\Div_{h_\cB}(\alpha).
$$
Indeed, since $(\cB, h_\cB)$ is assumed to be connected and without boundary, we obtain the following:

\begin{corollary}
	\label{Cor:Ricciflat}
	If $(\cB, h_\cB)$ is a compact Ricci-flat manifold then we always find a unique (up to a constant) $f_\cB \in C^\infty(\cB)$ s.\,t.\ the 
	Lorentzian manifold $(\cM^{(n + 2)}, g)$ as in Theorem \ref{Thm:Ricciflat} with $\cN = \cB \times \IS^1$ is Ricci-flat.
\end{corollary}

\begin{proof}	
	If $\omega = 0$ we choose $f_\cB$ such that $\alpha = -\tfrac{1}{4}df_\cB$ with $\alpha \in \omega$. Otherwise,
	since $\Div(\alpha) = \Div(\alpha^\sharp)$ and $\int_{\cB} \Div(\alpha^\sharp) = 0$ as $\partial \cB = \emptyset$, 
	we always find a unique (up to a constant) solution $f_\cB \in C^\infty(\cB)$ to the Poisson equation $\Delta^{h_\cB}(f_\cB) = -4\Div_{h_\cB}(\alpha)$,
	cf.\ \cite[Theorem 4.7]{1}.
\end{proof}

If $(\cB, h_\cB)$ is compact with $b_1(\cB) > 0$ and $\omega \neq 0$, the representative $\alpha \in \omega$ needs to be chosen non-harmonic for the function $f_\cB \in C^\infty(\cB)$ to be non-constant. Moreover, note that the condition $b_1(\cB) > 0$ is satisfied for compact Ricci-flat Riemannian manifolds whenever on $(\cB, h_\cB)$ exists at least one Killing vector field, since
in this case $b_1(\cB) = \dim_\R\LAkill(\cB, h_\cB)$, cf.\ \cite[Theorem 1.84]{4}. For example one may take the Ricci-flat metric on some Calabi-Yau manifold of dimension $n = 2m$, i.\,e.\ a compact Kähler manifold with trivial first Chern class. Examples are e.\,g.\ ${\rm K3} \times \T^k$ or more generally
products of compact Hyperkähler manifolds, i.\,e.\ a $4k$-dimensional Riemannian manifold with holonomy contained in $\operatorname{Sp}(k)$, with the flat torus. 
Another list of examples can be constructed from \cite[Theorem 4.1]{20}. \par
	Moreover, by Proposition \ref{Prop:CompleteInt}, the compact manifolds in Corollary \ref{Cor:Ricciflat} and thus in particular the just stated examples, are all complete.
	
\begin{corollary}
	\label{Cor:RicciflatComplete}	
	Every compact Ricci-flat Lorentzian manifold occurring in Corollary \ref{Cor:Ricciflat} is complete.
	This even holds for arbitrary $f_\cB \in C^\infty(\cB)$.
\end{corollary}

\begin{proof}
	Choose $\omega \in H_{\rm dR}^1(\cB)\cap H^1(\cB,\Z)$ and a representative $\alpha \in \omega$. If $\omega = 0$ then Proposition \ref{Prop:CompleteInt}\textit{(ii)} proves the statement.	
	When $\omega \neq 0$ and $(\cB, h_\cB)$ is assumed to be compact and Ricci-flat,
	we can write $\alpha = \widehat{\alpha} + d\varphi$, where $\widehat{\alpha} = K^\flat$ is the dual 1-form to a Killing field
	$$
		K \in \LAkill(\cB, h_\cB) = \{ X \in \Gamma(T\cM) \ | \ \nabla^{h_\cB}X = 0 \}.
	$$
	Let $\Psi := \alpha \wedge \eta$, $\widehat{\Psi} := \widehat{\alpha} \wedge \eta$
	and $A$, $\widehat{A}$ denote corresponding connection forms, i.e.\ with $dA = -2\pi i \pi^*\Psi$ and $d\widehat{A} = -2\pi i \pi^*\widehat{\Psi}$, respectively. Then
	\begin{equation}
		\label{proofRicciflatComplete1}
		A = \widehat{A} - 2\pi i (\varphi \circ \pi) \pi^*\eta.
	\end{equation}
	With the data chosen as in Theorem \ref{Thm:Ricciflat} we infer
	\begin{eqnarray*}
		g & = & 2i A \odot \pi^*\eta + (f + 1) \cdot \pi^*\eta \odot \pi^*\eta + \pi^*h_\cB \\
		  & \stackrel{(\ref{proofRicciflatComplete1})}{=} & 2i \widehat{A} \odot \pi^*\eta + (4\pi(\varphi \circ \pi) + f) \cdot \pi^*\eta \odot \pi^*\eta + \pi^*h_\cB \\
		  & = & 2i \widehat{A} \odot \pi^*\eta + \widehat{f} \cdot \pi^*\eta \odot \pi^*\eta + \pi^*h_\cB
	\end{eqnarray*}	
	for $\widehat{f} := 4\pi(\varphi \circ \pi) + f$. Hence, the Lorentzian manifold $(\cM^{(n + 2)}, g)$ is of type $(\widehat{\Psi}, \widehat{A}, \eta, \widehat{f})$ over $(\cN, h)$
	and	the assumptions of Proposition \ref{Prop:CompleteInt}\textit{(i)} are all satisfied, yielding the completeness.
\end{proof}
	
For the Einstein case with non-zero cosmological constant and particular Ricci-flat cases one has the following non-existence result:
	
\begin{proposition}
	\label{Prop:RicciflatSphereBundle}
	Let $(\cM^{(n + 2)}, g)$ be any Lorentzian manifold of type $(\Psi, A, \eta, f)$ over $(\cN, h)$. Then:	
	\begin{itemize}
		\item[(i)] $(\cM^{(n + 2)}, g)$ cannot be an Einstein manifold with non-zero cosmological constant.
		\item[(ii)] Let $(\cM^{(n + 2)}, g)$ be Ricci-flat and $\cN$ compact. If either
			\begin{itemize}
				\item[(a)] $\eta^\sharp$ is a $h$-Killing field, $\zeta(f) = 0$ and $\Psi(\eta^\sharp, \cdot) = 0$, or
				\item[(b)] $\eta^\sharp$ is $h$-parallel, $\cN = \cB \times \IS^1$ and $\Psi \in \Omega^2(\cB)$,				
			\end{itemize}			
			then $\Psi \in \Omega^2(\cN)$ must already vanish identically.
	\end{itemize}	 
\end{proposition}

\begin{proof}
	To prove \textit{(i)}, suppose the Lorentzian manifold $(\cM^{(n + 2)}, g)$ of type $(\Psi, A, \eta, f)$ over $(\cN, h)$ is an Einstein manifold. 
	Then, by Lemma \ref{lemma:Ric} (\ref{ref:equRic1}), the cosmological constant $\Lambda$ has to be equal to $\frac{1}{2}(\Hess f)(\xi, \xi)$.
	Hence, $(\Hess f)(\xi, \xi)$ has to be constant on each fiber since $0 = \xi(\Lambda) = \xi((\Hess f)(\xi, \xi))$	
	implying $(\Hess f)(\xi, \xi)|_{\pi^{-1}(y)} \equiv const$ for all $y \in \cN$. 
	As a consequence, such $f \in C^\infty(\cM)$ would give rise (by 
	passing to a local trivialization) to	a function $\hat f \in C^\infty(\IS^1)$ with constant Laplacian on $\IS^1$. Hence, $f$ is then contant
	on the fibers. But this is a contradiction to $\Lambda \neq 0$. \par
		To see \textit{(iia)} assume that $(\cM^{(n + 2)}, g)$ is Ricci-flat and the assumptions above hold true. Note that necessarily $\xi(f) = 0$.
	Then, by Proposition \ref{Prop:CompleteNonInt}, there exists a
	timelike Killing vector field $K \in \Gamma(T\cM)$. Due to \cite[Theorem 3.2]{28}, $K$ then has to be parallel. This is the case if and only if $\Psi$ vanishes, since
	$g(\nabla^g_{\e_i} K, \e_j) = \Psi_{ij}$ by Lemma \ref{equ:CovDer}c). \par
		Finally, in the case \textit{(iib)}, Ricci-flatness of $(\cM^{(n + 2)}, g)$ implies $\Delta_{h} f = 2||\psi||_h$ by Lemma \ref{lemma:Ric} (\ref{ref:equRic2}), where we regard
	$f$ as a function on $\cN$ which is feasible since $f$ is constant along the fibers. Since necessarily $\int_\cN \Delta_{h} f = 0$ we infer $||\psi||_h = 0$ and hence $\Psi = 0$.
\end{proof}

Note that this proposition implies in particular, that the toric type constructions in \cite{24} with compact base $\cN = \cB \times \IS^1$ and $(\cB, h_\cB)$ being Ricci-flat or Einstein cannot produce Ricci-flat or Einstein metrics on $\cM$ provided that $\Psi \in \Omega^2(\cB)$ is not chosen to be zero.
	
	\section{Holonomy}
\label{Section:Holonomy}

In \cite[Theorem 3]{7} there is given a criterion to compute the full holonomy group of a Lorentzian manifold $(\cM^{(n + 2)}, g)$ with parallel light-like vector field. We intend to apply this to the manifolds occurring in Theorem \ref{Thm:Ricciflat}. Therefore, assume we can show that the universal cover of the examples obtained from Theorem \ref{Thm:Ricciflat} is of the form $\widetilde{\cM} = \R^2 \times \cN$ while the metric on the universal cover is given as
$$
	\widetilde g_{(u,v,p)} = 2dudv + \kappa(u,p) du^2 + A_u \odot du + \Theta_p
$$
with $A = \{A_u\}$ a family of one-forms on $\cN$ and $\Theta$ a Riemannian metric on $\cN$. Although the 1-forms $A_u$ do not occur in \cite[Theorem 3]{7} it is -- by following the proof therein -- not hard to verify
that \textit{each isometry $\sigma$ of $(\R^2 \times \cN, \widetilde g)$ satisfies the assumptions} made in \cite[equation (13)]{7}, namely that
$$
	\sigma(u,v,p) = (a_\sigma^{-1}u + b_\sigma, a_\sigma v + \tau_\sigma(u,v,p), \nu_{\sigma}(u,v,p))
$$
with $a_\sigma \in \R^*, b_\sigma \in \R, \tau_\sigma \in C^\infty(\widetilde{\cM})$ with $\partial_v(\tau_\sigma) = 0$ and $\nu_\sigma : \widetilde{\cM} \To \cN$ such that 
$\partial_v(\nu_\sigma) = 0$ and $\nu(u,v,\cdot)$ is an isometry of $(\cN, \Theta)$ for all $u,v \in \R$. Then we can compute the full holonomy of $(\cM, g)$ by \cite[Theorem 3]{7} and obtain the following.

\begin{proposition}
	\label{prop:HolBLL}
	Under the assumptions above, it holds
	$$
		\Hol_x(\cM^{(n + 2)}, g) = Q \cdot \Hol_{\widetilde{x}}(\widetilde{\cM}^{(n + 2)}, \widetilde{g}) = Q \cdot \Hol_x^0(\cM^{(n + 2)}, g)
	$$
	where $\Phi : \widetilde{\cM} \To \cM$ denotes the universal covering, $\widetilde{x} = (u,v,p)$, $\Phi(\widetilde{x}) = x$, and
	$$
		Q := \left\langle Q(\sigma) \ | \ \sigma \in \pi_1(\cM) \right\rangle \subset \R^* \times {\rm O}(n)
	$$
	with $Q(\sigma) := (a_\sigma, d\mu_{\sigma^{-1}}^{-1} \circ \mathcal P_\sigma^\Theta)$. Here,	$\mu_\sigma := \nu_\sigma(u,v,\cdot)$	and 
	$\mathcal P_\sigma^\Theta$ the parallel transport w.r.t.\ $\Theta$ along some curve in $\cN$ from $p$ to $\mu_{\sigma^{-1}}(p)$.
\end{proposition}
	
Indeed, we find for the universal cover $(\widetilde{\cM}, \widetilde{g})$ of the Lorentzian manifolds of type $(\Psi, A, \eta, f)$ over $(\cN, h)$ appearing in Theorem \ref{Thm:Ricciflat} the following.

\begin{proposition}
	\label{Thm:univCover}
	Let $(\cM^{(n + 2)}, g)$ be a Lorentzian manifold of type $(\Psi, A, \eta, f)$ over $(\cN, h)$ for a smooth function $f \in C^\infty(\cM)$, constant along the fibers and
	$\cN = \cB \times \IS^1$ for a compact Riemannian manifold $(\cB, h_\cB)$ with $b_1(\cB) > 0$. Choose for $\eta$ the coordinate 1-form on $\IS^1$, $h = h_\cB \oplus \eta^2$ and 
	$\Psi := \alpha \wedge \eta$ for some nowhere vanishing closed 1-form $\alpha$ s.t.\ $[\alpha] \in H_{\rm dR}^1(\cB) \cap H^1(\cB,\Z)$. Then the universal cover $(\widetilde{\cM}, \widetilde{g})$ is 
	isometric to a manifold
	\begin{equation}
		\label{Equ:UnivMetric}
		(\R^2 \times \mathcal S, \ \Xi_{(u,v,p)} = 2dudv + \kappa(u,p)du^2 + A_u \odot du + \Theta_{p}),
	\end{equation}
	with $A_u = 2(u + a(s))ds$, $a \in \C^\infty(\R)$, and
	where $s$ is the $\R$-coordinate of $\mathcal S = \R \times \mathcal A$ which is the universal cover of a leaf of the integrable screen distribution $\IS|_{\mathcal L^\bot}$ defined in (\ref{equ:Z})
	on page \pageref{equ:Z}. Further, $\mathcal L^\bot$ is a leaf of $\IL^\bot$,
	$\kappa : \widetilde{\cM} \To \R$ is a smooth function not depending on the $v$-coordinate and $\Theta$ is a Riemannian metric on 
	$\mathcal S$ which coincides with the lift of $\pi^*h_\cB$ to the universal cover, restricted to $\mathcal S$.
\end{proposition}

\noindent Before we give the proof of the proposition, recall the following lemma \cite{21}.

\begin{lemma}
	\label{lem:splitting}
	Let $\cM$ be a manifold admitting a closed, nowhere vanishing one-form $\eta$. Assume that there is a complete vector field $Z$ such that $\eta(Z)=1$. 
	Then the leaves of the distribution $\ker(\eta)$ are all diffeomorphic to each other under the flow $\phi_t$ of $Z$, and the universal 
	cover $\widetilde{\cM}$ of $\cM$ is diffeomorphic to $\R \times {\cN}$ with the diffeomorphism given as
	$\R \times {\cN} \ni (u,p) \longmapsto \phi_u(p) \in \widetilde{\cM}$, where ${\cN}$ is the universal cover of a leaf of $\ker(\eta)$.
\end{lemma}

\begin{proof}[Proof of Proposition \ref{Thm:univCover}]
	To this end, let a tilde ahead of any object denote the lift to the universal cover. Moreover, we will use, locally, as a basis of $T\cB$ the $h_\cB$-orthonormal vector fields
	$E_\alpha, E_2, \ldots, E_n$ with $E_\alpha := \frac{\alpha^\sharp}{||\alpha^\sharp||}$ and $E_2, \ldots, E_n \in \ker\alpha$. 
	As usual we write $\e_i := E_i^*$ and write $\e_\alpha := E_\alpha^*$. \par
		We first show how to separate $\R^3$ from the universal cover $\widetilde{\cM}$ using Lemma \ref{lem:splitting}. 
	Indeed, $\pi^*\eta$ is closed on $\cM$ and $(\pi^*\eta)(\zeta) = 1$. Moreover,
	the 1-form $\pi^*\alpha$ on $\cM$ is closed, too, and fulfills $(\pi^*\alpha)(\e_\alpha) = 1$. Finally, fix a leaf 
	$\mathcal L^\bot$ of $\IL^\bot$. Since $\IS$ is horizontal and integrable, cf.\ Lemma \ref{lemma:Screen},
	the 1-form $\zeta^\flat = g(\zeta, \cdot) = \rmi A + (f + 1)\pi^*\eta$ is closed on $\mathcal L^\bot$. We can now apply Lemma \ref{lem:splitting} three times:
	$$
		\widetilde{\cM} 
		\stackrel{\pi^*\eta}{\simeq} \R \times \widetilde{\mathcal L^\bot} 
		\stackrel{\zeta^\flat}{\simeq} \R \times \R \times \mathcal S
		\stackrel{\pi^*\alpha}{\simeq} \R^3 \times \mathcal A,		
	$$
	where $\mathcal S$ is a fixed leaf of $\widetilde{\IS}|_{\widetilde{\mathcal L^\bot} }$ and $\mathcal A$ is a leaf of $\ker \widetilde{\pi^*\alpha}|_{\mathcal S}$. Recall that the diffeomorphisms
	are given by the flows of $\widetilde{\zeta}$, $-\widetilde{\xi}$ and $\widetilde{\e_\alpha}$, respectively. To be more precise, let
	$\{\varphi_u^\eta\}_{u \in \R}$, $\{\varphi_v^\xi\}_{v \in \R}$ and $\{\varphi_s^\alpha\}_{s \in \R}$ denote the corresponding flows of the latter vector fields, 
	respectively.	Then
	$$
		\Phi : \R^3 \times \mathcal A \ni (u,v,s,p) \longmapsto \varphi_u^\eta(\varphi_v^\xi(\varphi_s^\alpha(p))) \in \widetilde{\cM}
	$$
	is the asserted diffeomorphism. \par
		Since all vector fields except $\zeta$ and $\e_\alpha$ commute, we obtain 
	$$
		d\Phi(\partial_v) = -\widetilde{\xi}, \  d\Phi(\partial_u) = \widetilde{\zeta}.
	$$	
	As $\mathscr{L}_{\zeta}(\pi^*\alpha) = 0$, the flow of $\widetilde{\zeta}$ preserves $\widetilde{\pi^*\alpha}$ and thus
	$$
		\widetilde{g}(d\Phi(\partial_s), \widetilde{\e_\alpha}) 
		= \widetilde{g}(d\varphi_u^{\eta}(\widetilde{\e_\alpha}), \widetilde{\e_\alpha})
		= [(\varphi_u^{\eta})^*\widetilde{\pi^*\alpha}](\widetilde{\e_\alpha}) = 1.
	$$
	Moreover let, locally, $\omega_j := g(\e_j, \cdot)$, $j = 2, \ldots, n$. Then $\mathscr{L}_{\zeta}\omega_j = 0$ and hence 
	$$
		\widetilde{g}(d\Phi(\partial_s), \widetilde{\e_j})
		= \widetilde{g}(d\varphi_u^{\eta}(\widetilde{\e_\alpha}), \e_j)
		= (\varphi_u^{\eta})^*\widetilde{\omega_j}(\widetilde{\e_\alpha}) 
		= \widetilde{\omega_j}(\widetilde{\e_\alpha}) 
		= \omega_j(\e_\alpha) 
		= 0.
	$$	
	Therefore, we obtain that
	$$
		d\Phi(\partial_s) = \widetilde{\e_\alpha} + \tau \cdot \widetilde{\xi}
	$$
	for some $\tau \in C^\infty(\widetilde{\cM})$. Since $d\pi^*\alpha = 0$ we obtain $\mathscr{L}_{\e_\alpha}(\pi^*\alpha) = 0$. Hence, every flow defining $\Phi$
	preserves $\pi^*\alpha$ and since, locally, $\IS = \spann\{\e_2, \ldots, \e_n\}$, we see that
	$$
		d\Phi(\widetilde{\e_i}) \in \Gamma(\widetilde{\ker\pi^*\alpha}) \ i = 2, \ldots, n.
	$$	
	Since $\pi^*\Psi \in \Omega^2(\cM)$ is closed, its lift to the universal cover is exact. More precisely we have
	$$
		\Phi^*\widetilde{\pi^*\Psi} = \Phi^*\widetilde{\pi^*\alpha} \wedge \Phi^*\widetilde{\pi^*\eta} = ds \wedge du
	$$
	as $\Phi^*\widetilde{\pi^*\eta} = du$ and $\Phi^*\widetilde{\pi^*\alpha} = ds$. Hence,
	\begin{equation}
		\label{equ:prunivCover1}
		\rmi \Phi^*d\widetilde{A}
		= \rmi \Phi^*\widetilde{F^A}
		= 2 \Phi^*\widetilde{\pi^*\Psi} 
		= 2 ds \wedge du.
	\end{equation}
	Using this together with $\widetilde{\rmi A}(d\Phi(\partial_v)) = -\rmi A(\xi) = 1$ and 
	$\widetilde{\rmi A}(d\Phi(\partial_s)) = \tau \cdot \rmi A(\xi) = -\tau$, we see that
	$$
		\Phi^*(\widetilde{\rmi A}) = dv - \tau ds
	$$
	and hence $d\tau = -2du - b(s)ds$ by \eqref{equ:prunivCover1}, whence $\tau = -2(a(s) + u)$ for $2\tfrac{d}{ds}a = b$. Summarizing we get:
	\begin{eqnarray*}
		(\Phi^*\widetilde{g}) 
			& = & 2\Phi^*(\widetilde{\rmi A}) \odot \Phi^*\widetilde{\pi^*\eta} + (\widetilde{f} \circ \Phi + 1)(\Phi^*\widetilde{\pi^*\eta})^2 + \Phi^*(\widetilde{\pi^*h_\cB}) \\			
			& = & 2(dv + 2(u + a(s)) ds + (\widetilde{f} \circ \Phi + 1) du)du + \Phi^*(\widetilde{\pi^*h_\cB}) \\
			& = & 2dudv + \kappa du^2 + A_u \odot du + \Theta,
	\end{eqnarray*}
	where $\Theta := \Phi^*(\widetilde{\pi^*h_\cB})$ and $\kappa := \widetilde{f} \circ \Phi + 1$, while $\partial_v\kappa = 0$ 
	since $\xi(f) = 0$, i.\,e.\ $\kappa$ is independent of the $v$-coordinate.
\end{proof}

Next we need a description of the fundamental group of $\cM$ since this is contained in the groups $Q$ of Proposition \ref{prop:HolBLL}. Using Serre's long exact sequence
for the $\IS^1$-bundle $\pi : \cM \To \cN$ with $\cN = \cB \times \IS^1$ we obtain
\begin{equation}
	\label{Equ:Serre}
	\pi_2(\IS^1) = 0
	\rightarrow \pi_2(\cM)
	\stackrel{\varphi_1}{\To} \pi_2(\cB) = \pi_2(\cN)
	\stackrel{\varphi_2}{\To} \pi_1(\IS^1) = \Z
	\stackrel{\varphi_3}{\To} \pi_1(\cM)
	\stackrel{\varphi_4}{\To} \pi_1(\cB) \times \Z
	\rightarrow 0.
\end{equation}
This can be rewritten as the two short exact sequences
\begin{alignat}{7}
		0 \To \ & \pi_2(\cM) \ & \stackrel{\varphi_1}{\To} \ & \pi_2(\cB) & \stackrel{\varphi_2}{\To} \ & \im\varphi_2 \ & \To 0 
		\label{Equ:Serre-short} \\
		0 \To \ & \coker\varphi_2 \ & \stackrel{\varphi_3}{\To} \ & \pi_1(\cM) & \stackrel{\varphi_4}{\To} \ & \pi_1(\cB) \times \Z \ & \To 0
		\label{Equ:Serre-short-2}
\end{alignat}
To determine $\pi_1(\cM)$ from \eqref{Equ:Serre-short-2}, we make the following definition
\begin{definition}
	\label{def:split}
	We say that \dStr{$\pi_1(\cB)$ is split}, iff the short exact sequence \eqref{Equ:Serre-short-2} splits.
\end{definition}
\noindent For example, $\pi_1(\cB)$ is split, if it is a free group. We obtain:
\begin{proposition}
	\label{prop:fundamentalgroup}	
	If $\pi_1(\cB)$ is split then $\pi_1(\cM) \cong (\pi_1(\cB) \times \Z) \ltimes \coker\varphi_2 = (\pi_1(\cB) \times \Z) \ltimes \Z/\im\varphi_2$.
\end{proposition}

Since every subgroup of a free group is free, so is $\im\varphi_2 \subset \Z$ and consequently the sequence \eqref{Equ:Serre-short} always splits and
gives us a possibility to calculate either $\pi_2(\cB)$ or $\pi_2(\cM)$:

\begin{proposition}
	\label{prop:morealgebra}	
	$\pi_2(\cB) \cong \im\varphi_2 \ltimes \pi_2(\cM)$.
\end{proposition}

For example, in the easiest case where $\pi_2(\cB) = 0$ (e.g.\ when a cover of $\cB$ is contractible), then $\im\varphi_2 = 0$ and hence $\pi_1(\cM) = (\pi_1(\cB) \times \Z) \ltimes \Z$ by Proposition \ref{prop:fundamentalgroup}. If, for instance $\pi_2(\cB) = \Z$ (e.g.\ when $\cB = \C\PP^n$), then $\coker\varphi_2 \in \{1, \Z/k\Z, \Z\}$ and Propsition \ref{prop:morealgebra} may help to determine the correct case 
if one is able to get information about $\pi_2(\cM)$. For instance, if the leaves to $\IL^\bot = \xi^\bot$ are compact, then $\cM$ fibers over $\IS^1$ with each fiber diffeomorphic to a leaf $\mathcal L^\bot$ \cite[Corollary 8.6]{sharpe1997differential} and Serre's long exact sequence yields $\pi_2(\cM) \cong \pi_2(\mathcal L^\bot)$ and $\pi_1(\cM) = \Z \ltimes \pi_1(\mathcal L^\bot)$. \par
	We are now in the position to use Proposition \ref{prop:HolBLL} at the beginning of this section
to give a description of the holonomy of the Lorentzian manifolds of type $(\Psi, A, \eta, f)$ over $(\cN, h)$
considered in Proposition \ref{Thm:univCover}. We obtain:

\begin{theorem}
	\label{Prop:HolUniv}
	Let $(\cM^{(n + 2)}, g)$ be a Lorentzian manifold of type $(\Psi, A, \eta, f)$ over $(\cN, h)$ with the data chosen as in
	Proposition \ref{Thm:univCover} with $f \in C^\infty(\cN)$ s.t.\ $\Hess_B f|_\cB$ is non-degenerate in a point. Then the full holonomy group is given by	
	\begin{equation}
		\label{Eq:HolUniv}
		\Hol_x(\cM^{(n + 2)}, g) = O \cdot \Hol_q^0(\cB,h_\cB) \ltimes \R^n,
	\end{equation}	
	where $(\pr_\cB \circ \pi \circ \Phi)(\widetilde{x}) = q$, $\widetilde{x} = (u,v,p)$, $\Phi(\widetilde{x}) = x$ and 
	$$
		O := \left\langle (d\mu_{\sigma^{-1}})^{-1} \circ \mathcal P_\sigma^\Theta \ | \ \sigma \in \pi_1(\cM) \right\rangle \subset {\rm O}(n),
	$$
	with the notations as in Proposition \ref{prop:HolBLL}. Moreover, we can replace $\pi_1(\cM)$ by $\pi_1(\cB)$ in $O$, if $\pi_1(\cB)$ is split.
	In this case we	actually have
	\begin{equation}
		\label{Eq:HolUnivBsplit}
		\Hol_x(\cM^{(n + 2)}, g) = \Hol_q(\cB,h_\cB) \ltimes \R^n.
	\end{equation}
\end{theorem}

\begin{proof}
	The proof is threefold. As a first step we show that the manifolds occurring in Proposition \ref{Thm:univCover} have full holonomy
	$\Hol(\mathcal S, \Theta) \ltimes \R^n$ which is an easy adaption of the proof of \cite[Example 5.5]{6}. In a second step we prove
	that $\Hol(\mathcal S, \Theta)$ is isomorphic to $\Hol^0(\cB, h_\cB)$. Finally, we provide the arguments for the missing $\R^*$-factor in the
	groups $Q$ occurring in Proposition \ref{prop:HolBLL} and the fact that it suffices to consider generators $\sigma \in \pi_1(\cB)$. \par
		\textit{Step 1:} We prove that for the $(n + 2)$-dimensional manifold $\widetilde{\cM} = \R^2 \times \mathcal S$ equipped with the metric
	$\Xi_{(u,v,p)} = 2dudv + \kappa(u,p)du^2 + A_u \odot du + \Theta_{p}$	with simply-connected $\mathcal S \simeq \R \times \mathcal A$
	and $A_u = 2\mu ds := 2(u + a(s)) ds$, the full holonomy in the point
	$\widetilde{x} = (0,0,p)$ is given by
	\begin{equation}
		\label{equ:pHolUniv}
		\Hol_{\widetilde{x}}(\widetilde{\cM}, \Xi) = \Hol_{p}(\mathcal S, \Theta) \ltimes \R^n.
	\end{equation}
	Here, $p \in \mathcal S$ is a point s.t.\ $(\Hess_\Theta \kappa)(p)$ is non-degenerate.
	To prove \eqref{equ:pHolUniv} consider the basis $\partial_v,\partial_u,\partial_s, s_1,\ldots,s_{n-1}$ of $T\widetilde{\cM}$,
	where $\partial_s, s_1, \ldots, s_{n - 1}$ is a local $\Theta$-orthonormal frame in $T\mathcal S = \R \oplus \mathcal A$. Then, the only non-vanishing components of the Levi-Civita connection 
	$\nabla$ to $\Xi$ are given by
	\begin{equation}
		\label{equ:pHolUniv2}
		\begin{array}{lll}
			\nabla_{\partial_u}Y = \tfrac{1}{2}d\kappa(Y)\partial_v, & 
				\nabla_{\partial_u}\partial_u = (\tfrac{1}{2}d\kappa(\partial_u) + \tfrac{\mu}{2}d\kappa(\partial_s) - \mu)\partial_v + \partial_s -\tfrac{1}{2}\grad_\Theta\kappa, \\
			\nabla_{\partial_u}\partial_s = \tfrac{1}{2}d\kappa(\partial_s)\partial_v, &	\nabla_{\partial_s} X = \nabla_{X} \partial_s = \nabla^\Theta_{X} \partial_s, \\
			\nabla_{\partial_s}\partial_s = \nabla^\Theta_{\partial_s}\partial_s + a'(s)\partial_v, & \nabla_XY = \nabla_X^\Theta Y,			
		\end{array}
	\end{equation}
	where $X,Y \in \Gamma(T\mathcal A)$.
	Since the function $a$ does not depend on the $u$-coordinate we get for the curvature $R$ of $\Xi$
	\begin{equation}
		\label{equ:pHolUniv3}
		R(\partial_u,S_1)S_2 = -\frac{1}{2}\Hess_\Theta \kappa (S_1,S_2)\partial_v,		
	\end{equation}
	for all $S_1,S_2 \in \Gamma(T\cS)$. Hence, the holonomy algebra of $(\tem, \Xi)$ in $\widetilde{x}$ contains $\R^n$.
	Let $\gamma : [0,1] \To \tem$ be a curve with $\gamma(t) = (u(t), v(t), s(t), \delta(t))$
	with $\gamma(0) = (0,0,p)$, $p = (s,q)$ and $\delta : [0,1] \To \mathcal A$ a curve with $\delta(0) = q$.
	Then, for $X \in \Gamma(T\cS)$ being the $\Theta$-parallel vector field along $(s(t), \delta(t))$ with
	$X(0) = v \in T_p\cS$, we obtain for the parallel displacement $\mathcal P$ w.r.t.\ $\Xi$ that
	$$
		\mathcal P_{\gamma|_{[0,t)}}(v) = \varphi_v(t) \cdot (\partial_v \circ \gamma(t)) + X(t)
	$$
	with $\varphi_v : [0,1] \To \R$ defined as
	$$
		\varphi_v(t) = -\tfrac{1}{2} \int_0^t \left(\dot u(r) d\kappa_{\gamma(r)}(X(\gamma(r))) + \rho(r)\right) \rmd r,
	$$
	where
	$$
		\rho(r) = 
		\begin{cases}
			2\dot s(r) a'(s(r)), & v \in \R\partial_s, \\
			0, & v \in T_q\mathcal A.
		\end{cases}	
	$$
	Therefore, $\pr_{T_p\cS} \circ \mathcal P_\gamma|_{T_p\cS} = \mathcal P_{(s,\delta)}^\Theta$ which proves \eqref{equ:pHolUniv}. \par
	\textit{Step 2:} We are going to prove 
	\begin{equation}
		\label{equ:ScreenHolEqB}
		\Hol_x^0(\overline{\mathcal S}, \overline{h}) \cong \Hol_q^0(\cB, h_\cB).
	\end{equation}
	Here, $\overline{h} := \overline{\pi}^*h_\cB|_{\overline{\mathcal S} \times \overline{\mathcal S}}$, where $\overline{\mathcal S}$ is a leaf to the integrable screen distribution $\IS$ and
	$\overline{\pi}(x) = q$ with $x \in \overline{\mathcal S}$, where $\overline{\pi} : \overline{\mathcal S} \To \cB$ denotes the surjective map $\overline{\pi} := \pr_\cB \circ \pi |_{\overline{\mathcal S}}$.
	Then \eqref{equ:ScreenHolEqB} obviously implies $\Hol(\mathcal S, \Theta) \cong \Hol^0(\cB, h_\cB)$. \par
	First note that $(\overline{\mathcal S}, \overline{h})$ is geodesically complete since it is the restriction of a complete Riemannian metric $g^R$ on $\cM$
	(namely, $g^R = -A\odot A + \zeta^\flat \odot \zeta^\flat + \overline{\pi}^*h_\cB$) to a leaf (namely, $\overline{\mathcal S}$) of a foliation, cf.\ \cite[Exercise 10.4.28]{10}.	
	In addition it holds $\overline{\pi}^*h_\cB = \overline{h}$, i.e.\ $\overline{\pi}$ is a local isometry, and thus
	\begin{equation}
		\label{equ:PDisplLocIsom}
		d\overline{\pi}_x \circ \mathcal P^{\overline{h}}_{\overline{\gamma}} \circ d\overline{\pi}_x^{-1} = \mathcal P^{h_\cB}_{\overline{\pi} \circ \overline{\gamma}}
	\end{equation}
	for any loop $\overline{\gamma}$ in $x$. Finally, $\overline{\pi}$ is a Riemannian covering and hence every null-homotopic loop in $\cB$ lifts to a null-homotopic loop in $\overline{\mathcal S}$
	so \eqref{equ:ScreenHolEqB} follows from \eqref{equ:PDisplLocIsom}.	\par
	\textit{Step 3:} Let $\Phi : \widetilde{\cM} \simeq \R^2 \times \mathcal S \To \cM$ denote the universal covering from Proposition \ref{Thm:univCover} with $d\Phi(\partial_u) = \zeta$, $d\Phi(\partial_v) = -\xi$
	and hence $\Phi^*\pi^*\eta = du$. When $\sigma \in \pi_1(\cM)$ is a deck transformation of $(\widetilde{\cM}, \Xi = \Phi^*g)$, i.e.\ $\Phi \circ \sigma = \Phi$, then we see that
	$$
		\sigma^*du = \sigma^*(\Phi^*(\pi^*\eta)) = (\Phi \circ \sigma)^*\pi^*\eta = \Phi^*\pi^*\eta = du.
	$$
	Hence $u \circ \sigma = u + b_\sigma$, i.e.\ $a_\sigma = 1$ so there is no $\R^*$-factor in the groups $Q$ occurring in Proposition \ref{prop:HolBLL}. \par	
		Assume now that $\pi_1(\cB)$ is split, i.e.\ \eqref{Equ:Serre-short-2} splits, then 
	$\pi_1(\cM) \cong (\pi_1(\cB) \times \Z) \ltimes \Z/\im\varphi_2$ by Proposition \ref{prop:fundamentalgroup}. Let $x_0 \in \cM$. Then the integer factors in $\pi_1(\cM, x_0)$ come from the fundamental groups
	of the fibers and the circle in $\cN$. These are in turn generated by the flow of $\xi$ and $\zeta$ starting in $x_0$, respectively. Hence, if $\widetilde{x_0} = (u,v,p) \in \widetilde{M}$ with 
	$\Phi(\widetilde{x_0}) = x_0$ and for $k \in \Z$
	$$
		\widetilde{\gamma}_k^\zeta(t) := (u + kt, v, p), \ \widetilde{\gamma}_k^\xi(t) := (u, kt - v, p),
	$$
	then $\Phi \circ \widetilde{\gamma}_k^\zeta$ and $\Phi \circ \widetilde{\gamma}_k^\xi(t)$ are generators for the integer factors in $\pi_1(\cM, x_0)$ since it are integral curves
	of $k\zeta$ and $k\xi$, respectively. But neither $\widetilde{\gamma}_k^\zeta$ nor $\widetilde{\gamma}_k^\xi$ can connect $\widetilde{x_0}$ with $\sigma(\widetilde{x_0})$ for some
	isometry $\sigma$ of $(\widetilde{\cM}, \Xi)$ with $\nu(u,v,\cdot) \neq {\rm id}_{\mathcal S}$. So we can replace
	$\pi_1(\cM)$ by $\pi_1(\cB)$ in $O$. Since then	$\Hol_q(\cB,h_\cB) = O \cdot \Hol_q^0(\cB,h_\cB)$, cf.\ \cite[Proposition 3]{7}, this completes the proof.
\end{proof}

%

If $b_1(\cB) = 0$, we cannot choose a nowhere vanishing closed 1-form $\alpha \in \Omega^1(\cB)$ since it must be exact and hence $\alpha = d\tau$ for some smooth function $\tau \in C^\infty(\cB)$. But as $\cB$ was assumed to be compact, $\alpha = d\tau$ has at least one zero. Hence, Proposition \ref{Thm:univCover} cannot be applied in this case. However, if $b_1(\cB) = 0$, we may choose a different vector field
to split the first line from the universal covering. Indeed, if we choose the complete vector field $W := \zeta - 2(\tau \circ \pi)\xi$ on $\cM$ instead of $\zeta$, we can use the flow of its lift to the universal cover to split a line from $\widetilde{\cM}$ just as within the proof of Proposition \ref{Thm:univCover} but with the difference that now $[W, \e_\alpha] = 0$. We obtain the following.

\begin{proposition}
	\label{Thm:univCover2}
	Let $(\cM^{(n + 2)}, g)$ be a Lorentzian manifold of type $(\Psi, A, \eta, f)$ over $(\cN, h)$ as in Proposition \ref{Thm:univCover} but with $b_1(\cB) = 0$. 
	Then the universal cover $(\widetilde{\cM}, \widetilde{g})$ is isometric to a manifold
	\begin{equation}
		\label{Equ:UnivMetric2}
		(\R^2 \times \mathcal S, \ \Xi_{(u,v,p)} = 2dudv + \kappa(u,p)du^2 + \Theta_{p})
	\end{equation}
	with the notations as in Proposition \ref{Thm:univCover}.
\end{proposition}

\begin{proof}	
	Let $\alpha = d\tau$ and define $W := \zeta - 2(\tau \circ \pi)\xi$. Then $\pi^*\eta(W) = 1$ and the same methods as in the proof of 
	Proposition \ref{Thm:univCover} apply. Namely, by taking the flow $\{\varphi_u^W\}_{u \in \R}$ of $\widetilde{W}$ and $\{\varphi_v^\xi\}_{v \in \R}$ of $-\widetilde{\xi}$
	we can separate a line from $\widetilde{\cM}$ twice by Lemma \ref{lem:splitting}:
	$$
		\Phi : \R \times \R \times S \stackrel{\varphi_v^\xi}{\simeq} \R \times L^\bot \stackrel{\varphi_u^W}{\simeq} \widetilde{\cM}.
	$$
	Again, we will use, locally, as a basis of $T\cB$ the $h_\cB$-orthonormal vector fields	$E_\alpha, E_1, \ldots, E_n$ with $E_\alpha := \frac{\alpha^\sharp}{||\alpha^\sharp||}$ and $E_2, \ldots, E_n \in \ker\alpha$ and
	follow the notations in the proof of Proposition \ref{Thm:univCover}.	Then:
	$$
		[W, \e_\alpha] = \rmi F^A(\zeta,\e_\alpha)\xi + 2d\tau(\e_\alpha)\xi = -2 d\tau(E_\alpha)\xi + 2d\tau(E_\alpha)\xi = 0 \text{ and } [W, \e_i] = 0
	$$	
	for $i = 2, \ldots, n$. We obtain:
	$$
		d\Phi(\partial_u) = \widetilde{W}, \ d\Phi(\partial_v) = -\widetilde{\xi} \text{ and } d\Phi(\e_i) \in \Gamma(\widetilde{\ker\pi^*\eta}).
	$$
	The assertion now follows, since by the former equations,
	$$
		\Phi^*(\widetilde{iA}) = dv + 2(\widetilde{\tau \circ \pi} \circ \Phi) \cdot du.
	$$
	Setting $\kappa(u,p) := (\widetilde{f} + 1 + 4\widetilde{\tau \circ \pi}) \circ \Phi(u,0,p)$ completes the proof.
\end{proof}

Note that if not only $b_1(\cB) = 0$ but even $\cB$ is simply-connected, $\cM$ is diffeomorphic to $\T^2 \times \cB$ since
in this case the circle bundle is trivial as $[\Psi] = 0$. However, this must in general not be the case. 
Therefore it seems to be worthwhile to mention that the same conclusion about the holonomy as in Theorem \ref{Prop:HolUniv}
also holds for the case when $b_1(\cB) = 0$:

\begin{corollary}
	\label{Prop:HolUniv2}
	Let $(\cM^{(n + 2)}, g)$ be a Lorentzian manifold of type $(\Psi, A, \eta, f)$ over $(\cN, h)$ with the data chosen as in
	Proposition \ref{Thm:univCover2} with $f \in C^\infty(\cN)$ s.t.\ $\Hess_B f|_\cB$ is non-degenerate in a point. Then the full holonomy of
	$(\cM, g)$ is given as in Theorem \ref{Prop:HolUniv}.
\end{corollary}

So far we have just considered the case where the holonomy algebra of the Lorentzian manifolds
with special holonomy are of two certain types. Generally, the holonomy algebra of a Lorentzian
manifold with special holonomy, i.\,e.\ where its holonomy algebra acts indecomposable but 
non-irreducible, lies in the stabilizer of the invariant null line $L := W \cap W^\bot$
of the Lie algebra of $O(1,n+1)$, i.\,e.\
$$
	\hol_x(\cM^{(n + 2)},g) \subset \mathfrak{so}(1, n + 1)_L = (\R \oplus \mathfrak{so}(n)) \ltimes \R^{n}.
$$
It is well known \cite{34} that thus $\hol_x(\cM^{(n + 2)},g)$ can only be of four types.

\begin{theorem}
	\label{prop:holTypes}
	Let $\mathfrak{h} \subset \mathfrak{so}(1, n + 1)_L$ be an indecomposable subalgebra and let 
	$\mathfrak{g} := \pr_{\mathfrak{so}(n)}(\mathfrak{h})$ denote the orthogonal part. Then
	$\mathfrak{h}$ belongs to one of the following types:
	\begin{description}
		\item{Type 1:} $\mathfrak{h} = (\R \oplus \mathfrak{g}) \ltimes \R^{n}$,
		\item{Type 2:} $\mathfrak{h} = \mathfrak{g} \ltimes \R^{n}$,
		\item{Type 3:} $\mathfrak{h} = \{ (\varphi(X), X + Y, z) \ | \ X \in \mathfrak{z}(\mathfrak{g}), \ Y \in [\mathfrak{g},\mathfrak{g}], z \in \R^n \}$,
									 where $\varphi : \mathfrak{z}(\mathfrak{g}) \xtwoheadrightarrow{} \R$ is a surjective homomorphism,
		\item{Type 4:} $\mathfrak{h} = \{ (0, X + Y, \varphi(X) + z) \ | \ X \in \mathfrak{z}(\mathfrak{g}), \ Y \in [\mathfrak{g},\mathfrak{g}], z \in \R^k \}$,
									 where $\R^n = \R^m \oplus \R^k$, $0 < m < n$, $\mathfrak{g} \subset \mathfrak{so}(k)$ and
									 $\varphi : \mathfrak{z}(\mathfrak{g}) \xtwoheadrightarrow{} \R^m$ is a surjective homomorphism.
	\end{description}
\end{theorem}

Obviously, we have so far just considered Lorentzian manifolds of type $(\Psi, A, \eta, f)$ over $(\cN, h)$ which are of
type 1 or 2. By \cite[Proposition 6.2]{29}, they cannot be of type 3 since $R^{\nabla^\xi}(\rme_+,\xi) = 0$ if
and only if $\xi(f) = 0$. However, as we will see, for appropriate choices of the objects, we can obtain
Lorentzian manifolds of type 4 which are complete (but non-compact). We do not know if the other
existing examples \cite{31,32,35} for Lorentzian manifolds with holonomy of type 4 provide complete examples, too.
We use a characterization contained in \cite[Proposition 6.3]{29} which is a consequence of the Holonomy 
Theorem of Ambrose and Singer, and the curvature decomposition in \cite[Theorem 3.7]{33}.

\begin{proposition}
	\label{prop:BesvType4}
	A Lorentzian manifold $(\M^{(n+2)}, g)$ of type $(\Psi, A, \eta, f)$ over $(\cN, h)$ has type 4 holonomy algebra $\hol_x(\M, g)$ in $x \in \M$ if and only if
	there is a decomposition $\IS = \IS_1 \oplus \IS_2$ of a screen distribution such that:
	\begin{compactenum}[(i)]
		\item \label{prop:BesvType4-P1} $R^{\nabla^{\IS}}(X,Y)\Gamma(\IS_1) \subset \Gamma(\IS_1)$ and $R^{\nabla^{\IS}}(X,Y)\Gamma(\IS_2) = 0$
		for all $X,Y \in \Gamma(T\cM)$.
		\item \label{prop:BesvType4-P2} There exists a section $\varphi \in \Gamma(\operatorname{Hom}(\mathfrak{so}(\IS_1), \IS_2))$ s.t.\
		\begin{compactenum}[a)]
			\item \label{prop:BesvType4-P2a} $R^{\nabla^{\IS}}(X,Y) \in \ker\varphi$ for all $X,Y \in \Gamma(\IS)$,
			\item \label{prop:BesvType4-P2b} $\widehat{R}(\rme_+,\Gamma(\IS_2))\Gamma(\IS_2) = 0$
			and $\widehat{R}(\rme_+,X)Y = g(\varphi(R^{\nabla^{\IS}}(\rme_+,X)), Y)\xi$ for all vector fields $X \in \Gamma(\IS_1)$ and $ Y \in \Gamma(\IS_2)$,
			where $\widehat{R} = R^g - R^{\nabla^{\IS}}$.
			\item \label{prop:BesvType4-P2c} For any $y \in \cM$ and $\gamma : [0,1] \To \cM$ with $\gamma(0) = x$ and $\gamma(1) = y$ it holds 
			$$
				g_y( \varphi_y( R_y^{\nabla^{\IS}}(\rme_+, X) ), \mathcal P_\gamma^g(Y(x)) )
				 \! = \!
				  g_x( \varphi_x( \pr_{\IS_1} \circ \mathcal P_{\gamma^-}^g \circ R_y^{\nabla^{\IS}}(\rme_+, X) \circ \mathcal P_\gamma^g \circ \pr_{\IS_1} ), Y(x) )
			$$
			for arbitrary $X \in \Gamma(\IS_1)$ and $Y \in \Gamma(\IS_2)$.
		\end{compactenum}
	\end{compactenum}
\end{proposition}

Applying this to a certain family of Lorentzian manifold of type $(\Psi, A, \eta, f)$ over the flat manifold $\cN = \R^m \times \T^k$ gives us
the following.

\begin{proposition}
	\label{Ex:Type4}
	Let $(\cM^{(n + 2)}, g)$ be a Lorentzian manifold of type $(\Psi, A, \eta, f)$ over $(\cN, h)$, where we choose $\cN = \cB \times \IS^1$ with $\cB = \R^m \times \T^k$, 
	$\frac{k(k-1)}{2} \geq m > 0$ and $k \geq 2$. Denote by $\eta = du$ the coordinate 1-form on $\IS^1$ and fix a global trivialization of $T\cB$ by
	$\partial_1, \ldots, \partial_m,E_1,\ldots,E_k$. Furthermore, choose
	\begin{itemize}		
		\item $0 \neq [\Psi] \in H_{\rm dR}^2(\T^k) \cap H^2(\T^k,\Z)$ for a non-harmonic $\Psi$ and $[\Psi(x), \Psi(y)]_{\mathfrak{so}(k)} = 0$
		for all $x,y \in \T^k$, where $\Psi(x)$ is understood as an element of $\mathfrak{so}(k)$ w.r.t.\ the basis $E_1,\ldots,E_k$;
		
		\item smooth non-zero functions $\varphi_i : \R \To \R \setminus \{0\}$ with $\varphi_i(0) = 1$ for $i = 1,\ldots,m$;
		
		\item $h = h_\cB \oplus du$ with $h_\cB = \sum_{i = 1}^m \varphi_i^2 dx_i^2 \oplus h_{\T^k}$, where $h_{\T^k}$ is the flat metric on $\T^k$;
		
		\item $f := \widehat{f} \circ \pi$ for $\widehat{f} \in C^\infty(\cB)$ with 
		$\widehat{f}(y_1,\ldots,y_m,x) := -2 \sum_{(i,j) \in \Lambda} \Psi_{ij}(x)\Phi_{\lambda_i^j}(y_{\lambda_i^j})$, 
		where $\Phi_i$ is the antiderivative of $\varphi_i$ with $\Phi_i(0) = C_i \in \R$ and
		whereby $\lambda_i^j := \frac{(j-2)(j-1)}{2} + i$ and $\Lambda := \{(i,j) \in \{1,\ldots,k\}^2 \ | \ i < j, \ \lambda_i^j \leq m\}$.		
	\end{itemize}	
	Then $(\cM^{(n + 2)}, g)$ is a Lorentzian manifold with holonomy of type 4 and Abelian orthogonal part $\mathfrak{g} \subset \mathfrak{so}(k)$, where
	$\dim\mathfrak{g} = \left\lfloor \frac{\rank\psi}{2}\right\rfloor$.
\end{proposition}

\begin{proof}
	By the construction of $\cM$, we have that $\cM = \R^m \times \cM' \times \IS^1$ for the $\IS^1$-bundle $\pi : \cM' \To \T^k$ with $c_1(\cM') = [\Psi]$.	
	Let $\IS$ be the screen distribution	corresponding to the choice of the transversal vector field $Z \in \Gamma(T\cM)$ defined in (\ref{equ:Z})
	and $\nabla^{\IS} := \pr_\IS \nabla^g$ denote the induced connection from $\nabla^g$. By the choice
	of $\IS$ we have $\IS \simeq T\cB^* = \spann\{s_1^*, \ldots, s_n^*\}$ globally, where we set
	$s_i := \varphi_i^{-1}\partial_i$ and $s_j := E_j$ for $i = 1, \ldots, m$, $j = 1,\ldots, k$. Hence we have a splitting $\IS = \IS_2 \oplus \IS_1$ with
	$\IS_1 = (T\T^k)^*$ and $\IS_2 = \R^m$. \par	
		To this end we fix the point $x = (0, p, u) \in \cM$ for arbitrary $p \in \cM'$ and $u \in \IS^1$. Since the holonomy algebras in different points
	of the manifold are isomorphic, it suffices to prove that $\hol_x(\cM^{(n + 2)}, g)$ is of type 4. Computing $R^{\nabla^{\IS}}$ using Lemma~\ref{lemma:Curvature} we see that
	\begin{equation}
		\label{equ:Type4PrCurvS}
		R^{\nabla^{\IS}} = \nabla^h_{\cdot} \psi \wedge \pi^*\eta
	\end{equation}
	and thus
	$$
		R^{\nabla^{\IS}}(X,Y)\Gamma(\IS_1) \subset \Gamma(\IS_1), 
		R^{\nabla^{\IS}}(\Gamma(\IS),\Gamma(\IS)) = 0, 
		R^{\nabla^{g}}(X,Y)\Gamma(\IS_2) = R^{\nabla^{\IS}}(X,Y)\Gamma(\IS_2) = 0
	$$
	for all $X,Y \in \Gamma(T\cM)$. Therefore, in Proposition~\ref{prop:BesvType4} the properties \eqref{prop:BesvType4-P1}, \eqref{prop:BesvType4-P2a} and the first equation in \eqref{prop:BesvType4-P2b} are satisfied. \par
	We are left to choose a section $\varphi \in \Gamma(\operatorname{Hom}(\mathfrak{so}(\IS_1), \R^m))$
	for which \eqref{prop:BesvType4-P2b} and \eqref{prop:BesvType4-P2c} in Proposition~\ref{prop:BesvType4} hold. For every $y = (y_1,\ldots,y_m,q,v) \in \cM$,
	\begin{equation}
		\label{equ:Type4Pr1}
		\varphi_y : \mathfrak{so}(\mathcal (\IS_1)_y) \ni A_y \longmapsto \sum_{(i,j) \in \Lambda} A_{ij}(y)\varphi_{\lambda_i^j}(y_{\lambda_i^j}) \partial_{\lambda_i^j} \in \R^m
	\end{equation}
	defines a surjective linear map\footnote{
	For the purpose of clarifying this definition, we point out that the presented homomorphism $\varphi$ is nothing but the restriction of the canonical isomorphism $\mathfrak{so}(k) \cong \R^{k(k-1)/2}$
	given by the function $(a_{ij}) \mapsto (a_{12},a_{13},\ldots,a_{1k},a_{23},\ldots,a_{(k-1)k})$ to the first $m$ entries and weighted by the
	non-vanishing functions $\varphi_i$, $i = 1,\ldots, m$.
	}. To prove that Proposition~\ref{prop:BesvType4}~\eqref{prop:BesvType4-P2b} is satisfied we compute
	$$
		\widehat{R}(\rme_+,X)Y = \tfrac{1}{2}(\Hess_g f)(X,Y)\xi
	$$
	for all $X \in \Gamma(\IS_1)$ and $Y \in \Gamma(\IS_2)$. Moreover, we obtain for the Hessian of $f$
	\begin{equation}
		\label{equ:Type4Pr2}
		(\Hess_g f)_y(X, \partial_\ell) = X(\partial_\ell(f))(y)
		= -2X(\Psi_{i_0 j_0})\varphi_\ell(y_\ell)
		= -2g_y(\varphi(\nabla_{d\pi(X)}^h\psi), \partial_\ell)		
	\end{equation}	
	for $\lambda_{i_0}^{j_0} = \ell$ and all $X \in \Gamma(\IS_1)$ since $\varphi_\ell = \partial_\ell(\Phi_\ell)$. Therefore,		
	$$
		\widehat{R}(\rme_+,X)Y = g(\varphi(R^{\nabla^{\IS}}(\rme_+,X)), Y)\xi
	$$
	for all $X \in \Gamma(\IS_1)$ and $Y \in \Gamma(\IS_2)$ by \eqref{equ:Type4PrCurvS} and \eqref{equ:Type4Pr2} which proves Proposition~\ref{prop:BesvType4}~\eqref{prop:BesvType4-P2b}. \par
		Hence it remains to show that Proposition~\ref{prop:BesvType4}~\eqref{prop:BesvType4-P2c} holds. Let $\gamma : [0,1] \To \cM = \R^m \times \cM' \times \IS^1$ 
	be a path with $\gamma(0) = x$ and $\gamma(1) = y$ and write $\gamma(t) = (\delta^*(t), e^{\rmi u(t)})$. We define $\delta := \pi \circ \delta^*$.
	Furthermore, let $X(t) = a(t)\rme_+ + b(t)\xi(t) + Y(t)$ be a vector field along $\gamma$ with
	$Y \in \Gamma(\gamma^*(T\T^k\oplus\R^m))$. Then one computes for the parallel transport of any $E \in \{E_1^*, \ldots, E_k^*, \partial_1, \ldots, \partial_m\}$ 
	along $\gamma$ that
	$$
		\mathcal P^g_{\gamma}(E) = C_V \cdot \xi + V^*(1),
	$$
	where $C_V \in \R$ depends on $V \in \Gamma(T\T^k\oplus\R^m)$ which 
	is the solution to the ODE
	\begin{equation}
		\label{equ:Type4Pr3}
		\nabla_{\dot\delta}^h V = -\dot u \cdot \psi(V)
	\end{equation}
	with initial value $V(0) = d\pi(E)$. When $E = \partial_i$, then $\psi(E) = 0$ and
	we obtain the solution
	\begin{equation}
		\label{equ:Type4PrParallel1}
		V(t) = \mathcal P_{\delta|_{[0,t]}}^h(\partial_i(\delta(0))) = \frac{1}{\varphi_i(\delta_i(t))}\partial_i(\delta(t)).
	\end{equation}
	Hence, to solve \eqref{equ:Type4Pr3} we can write down \eqref{equ:Type4Pr3} as matrix equation of $(k \times k)$-matrices
	\begin{equation}
		\label{equ:Type4Pr5}
		\dot\Omega(t) = A(t) \cdot \Omega(t), \ \Omega(0) = \mathbb{I}_k
	\end{equation}
	for $A(t) := -\dot u(t) \cdot \psi(\delta(t))$, where $\mathbb{I}_k$ is the identity and $\psi$ is interpreted as an element of $\mathfrak{so}(k)$. We conclude that
	$\Omega(t) \in {\rm SO}(k)$ since $A(t) \in \mathfrak{so}(k)$. We obtain
	\begin{equation}
		\label{equ:Type4Pr6}
		\pr_{\IS_1} \, \circ \, \mathcal P_\gamma^g \circ \pr_{\IS_1} = \Omega_s(1) \in {\rm SO}(k)
	\end{equation}
	where $\Omega_s$ is the solution to \eqref{equ:Type4Pr5}. By \cite{30}, the solution $\Omega_s$ can explicitly
	written down as
	\begin{equation}
		\label{equ:Type4Pr7}
		\Omega_s(t) = \exp\left( \int_0^t A(\tau)\rmd\tau \right)
	\end{equation}
	since $[\Psi(\delta(\tau_1)), \Psi(\delta(\tau_2))]_{\mathfrak{so}(k)} = 0$ for all $\tau_1,\tau_2 \in [0,1]$ implying
	$[A(\tau_1), A(\tau_2)]_{\mathfrak{so}(k)} = 0$. Equation \eqref{equ:Type4Pr7} in turn implies that
	\begin{equation}
		\label{equ:Type4Pr8}
		\Omega_s(t_1)\Omega_s(t_2) = \Omega_s(t_2)\Omega_s(t_1)
	\end{equation}
	for $t_1,t_2 \in [0,1]$, which, by setting $t_2 = 1$ and differentiating in $t_1 = 1$, yields
	\begin{equation}
		\label{equ:Type4Pr9}
		\Psi(\delta(1))\Omega_s(1) = \Omega_s(1)\Psi(\delta(1)).
	\end{equation}	
	We are now in the position to prove Proposition~\ref{prop:BesvType4}~\eqref{prop:BesvType4-P2c}. 
	Consider the left hand side of the equation occurring in Proposition~\ref{prop:BesvType4}~\eqref{prop:BesvType4-P2c} for $X \in \Gamma(\IS_1)$, $Y = \partial_i$
	and $y = (y_1,\ldots,y_m,q,v) \in \cM$. We compute:
	\begin{eqnarray}
		g_y( \varphi_y( R_y^{\nabla^{\IS}}(\rme_+, X) ), \mathcal P_\gamma^g(\partial_i(x)) )
			& = &
			\varphi_i(y_i)^{-1} \cdot 
			g_y( \varphi_y( R_y^{\nabla^{\IS}}(\rme_+, X) ), \partial_i(y) ) \nonumber \\
			& = & 
			-\varphi_i(y_i)^{-1} \cdot 
			X(\Psi_{i_0j_0})(q)\varphi_i(y_i) \nonumber \\
			& = & -X(\Psi_{i_0j_0})(q) \label{equ:Type4Pr10}
	\end{eqnarray}
	where $i_0,j_0 \in \{1,\ldots,k\}$ such that $\lambda_{i_0}^{j_0} = i$. 
	To compute the right hand side,	define 
	$$
		B_{ij}(p) := g_p(\pr_{\IS_1} \circ \mathcal P_{\gamma^-}^g \circ R_y^{\nabla^{\IS}}(\rme_+, X) \circ \mathcal P_\gamma^g(E_i^*(p)), E_j^*(p)),
	$$
	such that it turns into
	$$
		g_x( \varphi_x( \pr_{\IS_1} \circ \mathcal P_{\gamma^-}^g \circ R_y^{\nabla^{\IS}}(\rme_+, X) \circ \mathcal P_\gamma^g \circ \pr_{\IS_1} ), \partial_i(x) ) = g_x(\varphi_x(B_x), \partial_i(x)).
	$$
	We compute
	\begin{eqnarray}
		B_{ij}(x) 
			& = & g_x(\mathcal P_{\gamma^-}^g \circ R_y^{\nabla^{\IS}}(\rme_+, X) \circ \mathcal P_\gamma^g(E_i^*(x)), E_j^*(x)) \nonumber \\
			& = & g_x(R_y^{\nabla^{\IS}}(\rme_+, X) \circ \mathcal P_\gamma^g(E_i^*(x)), \mathcal P_\gamma^g(E_j^*(x))) \nonumber \\
			& = & -g_y((\nabla_X^g\pi^*\Psi)(y) \circ \Omega_s(1)E_i^*(y), \Omega_s(1)E_j^*(y)) \nonumber \\
			& \stackrel{\eqref{equ:Type4Pr9}}{=} & -g_y(\Omega_s(1) \circ (\nabla_X^g\pi^*\Psi)(y)(E_i^*(y)), \Omega_s(1)E_j^*(y)) \nonumber \\
			& = & -g_y((\nabla_X^g\pi^*\Psi)(y)(E_i^*(y)), E_j^*(y)) \nonumber \\
			& = & -(\nabla_X^g\pi^*\Psi)_{ij}(y) \nonumber \allowdisplaybreaks[0] \\
			& = & -X(\Psi_{ij})(q). \label{equ:Type4Pr11} \allowdisplaybreaks[0]
	\end{eqnarray}
	Using this, we infer
	$$		
		g_x( \varphi_x(B_x), \partial_i(x) ) = B_{i_0j_0}(x) \varphi_i(0) \stackrel{\eqref{equ:Type4Pr11}}{=} -X(\Psi_{i_0j_0})(q).
	$$
	Taking into account \eqref{equ:Type4Pr10} this shows Proposition~\ref{prop:BesvType4}~\eqref{prop:BesvType4-P2c} and completes the proof of the proposition.
\end{proof}

\begin{remark}
	To our knowledge, up to now no compact examples of Lorentzian
	manifolds with holonomy algebra	of type 4 do exist.	
	Unfortunately we do not know, how to replace the $\R^m$ factor
	in $\cM$ by some compact manifold of dimension $m$ (e.g.\ the torus). The simplest
	idea is to try to choose periodic functions $\varphi_i$ such that
	their antiderivative is a periodic function. But since $\varphi_i$
	needs to be non-vanishing (i.e.\ either positive or negative),
	this is impossible.
\end{remark}

Under additional assumptions we get completeness of the latter manifolds producing examples
for geodesically complete Lorentzian manifolds with holonomy of type 4.

\begin{lemma}
	\label{Lem:CompleteType4}
	If the functions $\varphi_i$, $i = 1,\ldots, m$, and $[\Psi] \in H_{\rm dR}^2(\T^k) \cap H^2(\T^k,\Z)$ can be chosen, such that
	$(\R^m, \sum_{i = 1}^m \varphi_i^2dx_i^2)$ is complete and for each $u \in \R$ the solutions $s \mapsto \delta(s) \in \cB$ to the equation
	\begin{equation}
		\label{Equ:ExCompleteType4}
		\frac{\nabla^h \dot\delta}{ds}(s) = \frac{u^2}{2} \grad_h \widehat{f}(\delta(s)) - u \cdot \psi(\dot\delta)
	\end{equation}	
	are defined on the whole real line, then the Lorentzian manifold of type $(\Psi, A, \eta, f)$ over $(\cN, h)$
	in Proposition \ref{Ex:Type4} is complete.
\end{lemma}

\begin{proof}
	Let $\gamma : t \mapsto \gamma(t) \in \cM$ be a curve with $\gamma(t) = (\alpha(t), e^{\rmi u(t)})$, where
	$\alpha : t \mapsto \alpha(t) \in \R^m \times \cM'$ and define $\delta := \pi \circ \alpha$.
	Note that we write for $\pi$ the projection $\pi : \R^m \times \cM' \To \R^m \times \T^k$ to make notation short.
	Indeed, $\pi$ restricted to $\R^m$ is just the identity. We write 
	$$
		\dot\gamma(t) = \dot u(t) \partial_u + \dot\alpha(t) = \dot u(t) \partial_u + v(t)\xi(t) + dr_{\rho(t)}(\dot\delta^*(t))
	$$
	with $\delta^*$ denoting the horizontal lift of $\delta$ with $\delta^*(0) = \alpha(0)$ and $r_u : \cM' \To \cM'$ 
	the right action of $u \in \IS^1$ on $\cM'$, while	$\rho : \R \to \IS^1$ is defined through the equation $r_{\rho(t)}(\delta^*(t)) := \alpha(t)$. This yields
	\begin{equation}
		\label{Equ:ExCompleteType4Pr1}
		\tfrac{\nabla^g \dot\gamma}{dt}(t)
			= \ddot u(t) \partial_u + (\dot v(t) - \dot u(t)d\widehat{f}(\dot\delta))\xi(t) + \dot u(t)(\psi(\dot\delta) - \tfrac{1}{2}\dot u(t) \grad_h\widehat{f}) + \tfrac{\nabla^h \dot\delta}{dt}(t).
	\end{equation}
	Let $x = (y,p,e^{\rmi u_0}) \in \cM$ and $v \in T_x\cM$ be arbitrary with $v = u_1 \cdot \partial_u + \lambda \cdot \xi(x) + w$ where $w \in T_x\cB^* \cong \R^m \oplus \mathcal H_p\cM'$.\footnote{
	For any principal bundle $\mathcal P \To \cB$ we denote by $\mathcal H \mathcal P$ its horizontal bundle.}
	To prove completeness, we have to provide a geodesic $\gamma : \R \To \cM$ with $\gamma(0) = x$ and $\dot\gamma(0) = v$ defined on the whole
	real line. Let $\delta : \R \To \cB$ be a solution to \eqref{Equ:ExCompleteType4} for $u = u_1$ with $\delta(0) = \pi(y,p)$ and $\dot\delta(0) = d\pi(w)$.
	By \eqref{Equ:ExCompleteType4Pr1}, for each geodesic $\gamma$ we have $u(t) = u_1 t + u_0$. Therefore, 
	$\tau(t) := \dot u(t)d\widehat{f}(\dot\delta)$ is defined on the whole $\R$ and we define by $\mathcal T : \R \To \R$ its antiderivative with
	$\mathcal T(0) = \lambda$. If $T : \R \To \R$ is the antiderivative of $\mathcal T$ with $T(0) = 0$, then we define by
	$$
		\alpha(t) := r_{\rho(t)}(\delta^*(t))
	$$
	for $\rho(t) := e^{\rmi T(t)}$ a curve in $\R^m \times \cM'$ with $\delta^*$ denoting the horizontal lift of $\delta$ with $\delta^*(0) = (y,p)$. 
	We do now claim	that 
	$$
		\gamma(t) := (\alpha(t), e^{\rmi (u_1 t + u_0)})
	$$
	is the required geodesic with $\gamma(0) = x$ and $\dot\gamma(0) = v$. To see this, first note that
	$$
		\gamma(0) = (\alpha(0), e^{\rmi u_0}) = (\delta^*(0), e^{\rmi u_0}) = x.
	$$
	In order to verify that $\gamma$ is a $g$-geodesic, recall \eqref{Equ:ExCompleteType4Pr1} and the formula
	$$
		\frac{d}{dt} \alpha(t) = dr_{\rho(t)}(\tfrac{d}{dt}\delta^*(t)) + X(\alpha(t)),
	$$
	where $X \in \Gamma(T\cM')$ is the fundamental vector field to $dL_{\rho(t)^{-1}}(\dot\rho(t)) \in \rmi\R$ with
	$L_u : \IS^1 \To \IS^1$ denoting the left-multiplication by $u$ in $\IS^1$. In fact, 
	$$
		X(\alpha(t)) = \mathcal T(t) \cdot \xi(t),
	$$
	while $\pi \circ \alpha = \delta$. Hence $\dot v(t) = \tau(t) = \dot u(t)d\widehat{f}(\dot\delta)$
	and since $\delta$ satisfies \eqref{Equ:ExCompleteType4}, $\gamma$ is a $g$-geodesic with $\dot\gamma(0) = v$.	
\end{proof}

The following result provides an example for the existence of the required functions $\varphi_i$ and the 2-form $\Psi$ such that \eqref{Equ:ExCompleteType4} in Lemma \ref{Lem:CompleteType4}
is satisfied.

\begin{proposition}
	\label{Ex:CompleteType4}
	Let $\varphi_i \equiv 1$, $i = 1,\ldots, m$, and for $l = \left\lfloor \frac{k}{2} \right\rfloor$ define
	\begingroup
	\small
	$$
		\Psi(x_1,\ldots,x_k) := 
		\begin{pmatrix}
			0 						   & \chi_1(x_1,x_2) & \ldots & 0 							 & 0								\\
			-\chi_1(x_1,x_2) & 0						   &				&   							 & \vdots					  \\
			\vdots				   & 							   & \ddots	&  								 & 0								\\
			0							   &  						   & 				& 0 							 & \chi_l(x_{2l-1},x_{2l}) 	\\
			0							   &  						   & 				& -\chi_l(x_{2l-1},x_{2l}) & 0 								\\
		\end{pmatrix}
		\in \mathfrak{so}(2l)
	$$
	\endgroup
	for periodic functions $\chi_i : \T^2 \To \R$, $i = 1, \ldots, 2l$ such that $[\Psi] \in H_{\rm dR}^2(\T^k) \cap H^2(\T^k,\Z)$.\footnote{
	Note that the constructed $\Psi$ is of the form $\Psi = \Psi_1 + \ldots + \Psi_l$ where each $\Psi_i \in \Omega^2(\T^2)$
	is simply a 2-form on $\T^2$ not depending on the other coordinates. Hence we obtain $[\Psi] \in H^2(\cB,\Z)$ iff the integral over the fundamental class $[\T^2]$ for each $\Psi_i$ is an integer. 
	For example one may choose $\chi_i(x_i, x_{i+1}) := \sin(x_i)\sin(x_{i+1})$.}
	Then the Lorentzian manifolds of type $(\Psi, A, \eta, f)$ over $(\cN, h)$ provided in Proposition \ref{Ex:Type4} are geodesically complete.
\end{proposition}

\begin{proof}
	Let $\delta(s) = (\alpha_1(s), \ldots, \alpha_m(s),\beta_1(s), \ldots, \beta_k(s))$  be a path in $\cB = \R^m \times \T^k$
	and set $\alpha(s) = (\alpha_1(s), \ldots, \alpha_m(s))$ and $\beta(s) = (\beta_1(s), \ldots, \beta_k(s))$.	Since $h_\cB = h_1 \oplus h_{2}$ with $h_1 = \sum_{i = 1}^m \varphi_i^2dx_i^2$
	and $h_2 := h_{\T^k}$, equation	\eqref{Equ:ExCompleteType4} becomes
	\begin{equation}
		\label{Ex:CompleteType4Pr1}
		\left.\begin{array}{rcl}
			\ddot\alpha(s) & = & \tfrac{u^2}{2} \grad_{h_1} \widehat{f}(\delta(s)), \\
			\frac{\nabla^{h_2} \dot\beta}{ds}(s) & = & \tfrac{u^2}{2} \grad_{h_2} \widehat{f}(\delta(s)) - u\psi(\dot\beta(s)).
		\end{array}\right\}
	\end{equation}
	Taking into account the definition of $\widehat{f} \in C^\infty(\cB)$ in Proposition \ref{Ex:Type4} and $\Phi_i(x) = x + C_i$
	for $i = 1,\ldots, m$, equation \eqref{Ex:CompleteType4Pr1} turns into
	\begin{equation}
		\label{Ex:CompleteType4Pr2}
		\left.\begin{array}{rcl}
			\ddot\alpha_a(s) & = & -u^2 \Psi_{i_0j_0}(\beta(s)), \ a = 1,\dots,m, \\
			\frac{\nabla^{h_2} \dot\beta}{ds}(s) & = & 
			-u^2 \sum_{b = 1}^k \sum_{(i,j) \in \Lambda} \{E_b(\Psi_{ij})(\beta(s)) \cdot (\alpha_{\lambda_i^j}(s) + C_{\lambda_i^j})\}E_b 
			- u\psi(\dot\beta(s)).
		\end{array}\right\}
	\end{equation}
	where $i_0,j_0 \in \{1,\ldots,k\}$ such that $\lambda_{i_0}^{j_0} = a$. By integrating the first equation of \eqref{Ex:CompleteType4Pr2}
	twice and substituting this into the second equation we obtain equivalently:
	\begin{equation}
		\label{Ex:CompleteType4Pr3}
		\frac{\nabla^{h_2} \dot\beta}{ds}(s) =
		-u^4 \sum_{b = 1}^k \sum_{(i,j) \in \Lambda} \left\{E_b(\Psi_{ij})(\beta(s)) \left(C_{\lambda_i^j} - \int_0^s\!\int_0^t \Psi_{ij}(\beta(\tau)){\rm d}\tau{\rm d}t\right)\right\}E_b 
			- u\psi(\dot\beta(s))
	\end{equation}	
	By lifting this equation to $\R^k$, we obtain a second order non-linear differential equation of the form $y''(s) = F(s,y,y') := A(s, y(s)) + B(y(s))y'(s)$. 
	Since the partial derivatives of $A : [a,b] \times \R^k \To \R$ and $B : \R^k \To \R$ are bounded, $F : \R^{2k} \To \R^k$ is globally Lipschitz continuous and 
	\eqref{Ex:CompleteType4Pr3} exhibits a global solution $\widetilde{\beta} : \R \To \R^k$. Taking $\beta := \pi \circ \widetilde{\beta}$	for the canonical 
	projection $\pi : \R^k \To \T^k$ then yields the global solution on the torus.
\end{proof}

\noindent Combining Proposition \ref{Ex:Type4} and Proposition \ref{Ex:CompleteType4} we finally obtain the following result.

\begin{theorem}
	\label{Thm:Type4CplExamples}
	For each Abelian Lie subalgebra $\mathfrak{g} \subset \mathfrak{so}(k)$ there exists a complete indecomposable Lorentzian manifold with holonomy of type 4
	possessing $\mathfrak{g}$ as orthogonal part.
\end{theorem}

	\small

\bibliographystyle{amsalpha}
\bibliography{Bibliography}

\vspace{1cm}

\noindent\textsc{Daniel Schliebner}\newline
Humboldt-Universität zu Berlin, Institut für Mathematik\newline
Rudower Chaussee 25, Room 1.304, D--12489 Berlin, Germany.\newline
E-Mail: {\tt schliebn@mathematik.hu-berlin.de}

\end{document}